\documentclass[onecolumn, journal]{IEEEtran} % 单栏核心参数
\usepackage{amsmath,amsfonts}
\usepackage{algorithmic}
\usepackage{array}
\usepackage[caption=false,font=normalsize,labelfont=sf,textfont=sf]{subfig}
\usepackage{textcomp}
\usepackage{stfloats}
\usepackage{url}
\usepackage{float}
\usepackage{verbatim}
\usepackage{amsthm}  
\newtheorem{proposition}{Proposition}
\usepackage{subcaption} % For subfigures
\usepackage{mathtools}
\usepackage{graphicx}% 图片加载宏包
\usepackage{tabularx}   % 自适应宽度表格
\usepackage{multirow}   % 多行/列表格
\usepackage{hyperref} % 超链接宏包
\hypersetup{colorlinks=true, citecolor=blue} % 配置超链接颜色
\usepackage{cite}
\usepackage{times}
\usepackage{xcolor}
\usepackage{amssymb}
\usepackage{booktabs}
\usepackage{epstopdf}
\newtheorem{theorem}{Theorem}
\def\BibTeX{{\rm B\kern-.05em{\sc i\kern-.025em b}\kern-.08em
		T\kern-.1667em\lower.7ex\hbox{E}\kern-.125emX}}
\usepackage{balance}
\usepackage{fancyhdr}
\fancypagestyle{firstpage}{
	\fancyhf{} % 清空所有页眉页脚
	\fancyfoot[R]{\thepage} % 首页页码在右上角
}
\newtheorem{remark}{Remark}
\newtheorem{lemma}{Lemma}
\usepackage{etoolbox}
\makeatletter
\patchcmd{\IEEEbiography}
{\addvspace{0.5\baselineskip}}
{\vspace{-3mm}}
{}{}
\patchcmd{\endIEEEbiography}
{\addvspace{0.5\baselineskip}}
{\vspace{-3mm}}
{}{}
\makeatother

\DeclareGraphicsExtensions{.eps,.pdf,.png,.jpg}

\begin{document}
	\thispagestyle{firstpage}
	\title{Randomized Optimal Control Problem via Reinforcement Learning under Knightian Uncertainty\thanks{This work was supported by the National Natural Science Foundation of China (Grant Nos. 72301173, 62273003). (\textit{Corresponding author: Chen Fei.})}}
	\author{
		Ziyu Li, 
		Chen Fei, 
		Weiyin Fei
		\thanks{Ziyu Li and Weiyin Fei are with the School of Mathematics-Physics and Finance, Anhui Polytechnic University, Wuhu, Anhui 241000, China (e-mail: 1790941318@qq.com; wyfei@ahpu.edu.cn).}
		\thanks{Chen Fei is with the Business School, University of Shanghai for Science and Technology, Shanghai 200093, China (e-mail: chenfei@usst.edu.cn).}
	}
	
	\maketitle % 标题生成
	
	\indent\textbf{Abstract} % 顶格加粗显示abstract
		Considering that the decision-making environment faced by reinforcement learning (RL) agents is full of Knightian uncertainty, we describe the exploratory state dynamics equation to study the entropy-regularized relaxed stochastic control problem in a Knightian uncertainty environment. By employing stochastic analysis theory and the dynamic programming principle under nonlinear expectation, we derive the Hamilton-Jacobi-Bellman (HJB) equation and solve for the optimal policy that achieves a trade-off between exploration and exploitation. Subsequently, for the linear-quadratic (LQ) case, we prove that the optimal randomized feedback control follows a Gaussian distribution. Furthermore, we investigate how the degree of Knightian uncertainty affects the variance of the optimal feedback policy. Finally, we establish the solvability equivalence between non-exploratory and exploratory LQ problems under Knightian uncertainty and analyze the associated exploration cost.


	\section{Introduction}
	\IEEEPARstart{O}{ptimal} control has been widely applied in various fields, including engineering, economics, and biomedical sciences. Traditional optimal control methods, such as dynamic programming and Pontryagin's maximum principle, while possessing a rigorous mathematical foundation in theory, often encounter challenges such as high computational complexity or difficulties in modeling when dealing with high-dimensional state spaces, nonlinear systems, or environmental uncertainties. As a data-driven adaptive optimization method, reinforcement learning (RL) has recently achieved breakthrough progress in domains such as compliant control~\cite{Li2024}, multi-agent systems~\cite{Liu2024, Dai2025}, smart grid \cite{Dai2020a}, intelligent manufacturing \cite {qin2024value} and cybersecurity \cite{Dai2020b}, providing novel approaches for solving complex optimal control problems. RL learns optimal policy through agent-environment interactions without requiring precise system dynamics models, making it particularly suitable for scenarios where models are unknown or difficult to formulate and demonstrates strong learning capabilities in high-dimensional and nonlinear systems~\cite{jiao2024novel}.
	
	Despite the rapid development and wide application of RL techniques in solving optimal control problems, existing research still suffers from two major limitations. First, the majority of work has focused on RL algorithms for discrete-time markov decision processes, while paying considerably less attention to problems in continuous time and space. In practical scenarios such as high-frequency trading and autonomous driving, agents are required to interact with the environment at ultra-high frequencies or in a continuous manner, a significant drawback of addressing these problems in a discrete-time setting is that the resulting algorithms are highly sensitive to the time step size 
	\cite{yildiz2021continuous}. Second, most existing studies on continuous-time RL are confined to deterministic systems \cite{lee2021policy}, or although having introduced stochastic diffusion models, still lack data-driven solutions. Even parameter estimation methods remain vulnerable to model specification errors \cite{basei2022logarithmic}.
	
	Systematic exploration of RL for continuous-time diffusion processes, particularly data- or sample-driven solutions, has only recently gained momentum. Wang et al. \cite{wang2020reinforcement} proposed an entropy-regularized relaxed stochastic control framework for studying RL in continuous time and space, deriving the Boltzmann distribution as the general optimal randomized policy for environment exploration and action generation. Jia and Zhou \cite{jia2022policya} addressed the policy evaluation problem---learning the value function of a given stochastic policy---by formulating it as a martingale problem. Subsequently, Jia and Zhou \cite{jia2022policyb} investigated the policy gradient problem, expressing the gradient of the value function with respect to a parameterized stochastic policy as the expected integral of an auxiliary running reward function. Further extending this line of research, Jia and Zhou \cite{jia2023q} examined the $q$-learning aspect of continuous-time RL. Wang and Zhou \cite{wang2020continuous} applied the entropy-regularized relaxed stochastic control framework to mean-variance (MV) portfolio selection problems with linear-quadratic (LQ) structures, designing algorithms for extensive simulation and empirical studies. Dai et al. \cite{dai2023learning} advanced this research by investigating equilibrium MV policies that resolve the time-inconsistency in MV portfolio problems. Guo et al. \cite{guo2022entropy} extended Wang et al. \cite{wang2020reinforcement}'s framework to mean-field games, while Gao et al. \cite{gao2022state} adapted their methodology to the non-learning context of simulated annealing, modeling temperature control of Langevin diffusion for non-convex optimization.
	
	The theory of RL for continuous-time stochastic control has achieved remarkable progress since the foundational work by Wang et al. \cite{wang2020reinforcement}. These studies primarily focus on regular control problems within (jump-) diffusion processes, employing the Feynman-Kac formula and partial differential equations of value functions to derive gradient representations of value functions with respect to parameterized policies. However, Wang et al. \cite{wang2020reinforcement} and subsequent extensions did not consider the uncertainty of the model itself, namely, Knightian uncertainty. Knight \cite{knight1921risk} emphasized that the inherent uncertainty in probabilistic statistical models is fundamental and irreducible. Consequently, practical applications inevitably require consideration of Knightian uncertainty effects, and RL research is no exception. Generally speaking, the concept of linear expectation operates within the axiomatic framework of probability theory established by the eminent mathematician Kolmogorov in the 1930s. However, Peng \cite{peng2010nonlinear} demonstrated that problems exhibit increasingly dynamic characteristics and uncertainties as societal complexity grows, necessitating more fundamental theories and robust analytical tools. While classical probability theory effectively handles scenarios where models can be precisely determined through statistical methods, it fails when facing unknown probabilistic models. In contrast, nonlinear expectations enable robust quantitative analysis and computation of risks even when the underlying probability model remains undetermined.
	
	Currently, numerous scholars primarily focus on financial and economic issues under Knightian uncertainty. Chen and Epstein \cite{chen2002ambiguity} analyzed the distinction between risk premiums and ambiguity premiums using a continuous-time intertemporal model with multiple priors, establishing an asset pricing theory within the Knightian uncertainty framework. Fei et al. \cite{fei2019dynamics} investigated continuous-time principal-agent problems under Knightian uncertainty. Fei et al. \cite{fei2020contract} further integrated Knightian uncertainty with adverse selection to examine moral hazard in continuous-time principal-agent relationships. Ling et al. \cite{ling2022robust} proposed a dynamic robust contracting model to capture the principal’s ambiguity aversion toward the agent’s information disclosure, deriving optimal dynamic incentive-compatible contracts and endogenous corporate dividend policies. Niu and Zou \cite{niu2024robust} constructed a continuous-time dynamic multi-agent contracting model, demonstrating an inverted U-shaped relationship between overdetermination degree and team size. Sun et al. \cite{sun2023impact} explored how ambiguity preferences influence market participation and asset pricing. Building upon the properties of $G$-Brownian motion in sublinear expectation spaces, Fei \cite{chen2021optimal} developed an optimality principle for stochastic control problems under the nonlinear expectation framework. 
	
	The aforementioned studies systematically reveal the impact of Knightian uncertainty on financial-economic decision-making, demonstrating the limitations of the traditional probabilistic framework across domains ranging from asset pricing to contract theory. Notably, these limitations manifest more prominently in RL contexts: RL agents typically operate in complex, dynamically evolving environments replete with unquantifiable uncertainties and unknown factors. Taking autonomous vehicles as an example, uncertainties in road conditions, traffic flow, and weather patterns, some of which cannot be fully captured through historical data or statistical methods, necessitate explicit consideration of Knightian uncertainty. Neglecting such Knightian uncertainty factors, particularly the agent's potential decision-making behaviors under worst/best-case scenarios, will inevitably compromise the accuracy of environmental state predictions and reward estimations, thereby undermining both policy optimization effectiveness and decision-making robustness. Consequently, when investigating RL-based optimal control problems, it becomes imperative to appropriately incorporate Knightian uncertainty factors to bridge the sim-to-real gap. 
	
	Based on the aforementioned research literature, this paper considers extending the entropy-regularized relaxed stochastic control framework proposed by Wang et al. \cite{wang2020reinforcement} to the Knightian uncertainty environment to study the RL-based optimal control problem under Knightian uncertainty. This brings about several key theoretical challenges in the following aspects:
	
	\begin{enumerate}
		\item The framework of Wang et al. \cite{wang2020reinforcement} relies on linear expectation based on a single fixed probability measure, whereas the essence of Knightian uncertainty lies in the existence of a family of possible probability measures. Consequently, it is not feasible to study the RL-based optimal control problem under Knightian uncertainty within the linear expectation framework.
		
		\item Under sublinear expectation, the classical strong law of large numbers is no longer applicable. Therefore, the exploratory state dynamics cannot be derived by simply applying the classical strong law of large numbers.
		
		\item In the linear expectation framework, the corresponding HJB equation can be derived using the dynamic programming principle and Itô's formula. In contrast, under sublinear expectation, the resulting HJB equation will incorporate a nonlinear operator, and both its derivation logic and mathematical tools differ fundamentally from the classical case.
		
		\item Under sublinear expectation, it is not possible to prove the solvability equivalence between the exploratory and non-exploratory problems by comparing the analytical solutions of ordinary differential equations (ODEs).
	\end{enumerate}
	
	To address these challenges, this paper constructs a unified analytical framework based on sublinear expectation and achieves the following core theoretical results:
	
	\begin{itemize}
		\item We introduce the sublinear expectation theory and \textit{G}-Brownian motion proposed by Peng \cite{peng2010nonlinear} to characterize the uncertainty inherent in the probability model itself.
		
		\item We rigorously derive the exploratory state dynamics using the strong law of large numbers for sublinear expectations (Chen \cite{Chen2016}). Furthermore, we formulate the $G$-HJB equation—which incorporates a nonlinear operator $\widetilde{G}[\cdot]$ dependent on the uncertainty volatility interval—by applying the nonlinear dynamic programming principle \cite{chen2021optimal} and the $G$-Itô's formula.
		
		\item For the LQ case, we prove that the variance of the optimal policy explicitly depends on the Knightian uncertainty parameters, thereby revealing the novel phenomenon of ``adaptive exploration intensity varying with the level of uncertainty.''
		
		\item Within the sublinear expectation framework, we prove the solvability equivalence between the exploratory and non-exploratory problems by employing the $G$-Itô's formula and Gronwall's inequality.
		
		\item We demonstrate that the exploration cost is exclusively determined by the agent's intrinsic parameters and is independent of environmental uncertainty, and we provide convergence analysis for the policy.
	\end{itemize}

	Classical risk-sensitive control \cite{anugu2025ergodic}, $H_\infty$ control \cite{ma2024open}, modern distributionally robust RL (DRRL) \cite{shi2024distributionally}, and the $G$-expectation framework employed in this paper are all designed to address control under uncertainty. However, they operate in distinct environments and utilize fundamentally different methodologies. Classical risk-sensitive control and $H_\infty$ control both assume the existence of a true and known probabilistic model, essentially performing optimization within a known ``risk'' landscape without addressing the uncertainty inherent in the model itself.
	
	Both reference \cite{shi2024distributionally} and our work solve RL problem under the model uncertainty. Specifically, the significant differences between reference \cite{shi2024distributionally} and our work are as follows. First, our work is based on continuous time context, where we can employ the continuous time stochastic analysis such as It\^o's calculus, stochastic differential equations (SDE), nonlinear dynamic programming principle etc. However, reference \cite{shi2024distributionally} adopts the framework of discrete time, states and actions. This leads to the significant different methods of study. Second, the construct of uncertainty family is different, where the prior set of uncertainty probability measures in our work is closely related to $G$-Brownian motion with the degree of volatility uncertainty which can be estimated by employing $\varphi$-max-mean algorithm put forward by Peng \cite{peng2017}. The set of uncertainty family in reference \cite{shi2024distributionally} is based on a reference (nominal) probability measure, where the distance between other probability measure and reference probability measure is defined by Kullback-Leibler divergence. Third, our work employs the entropy regularized factor as a (reward) penalized term while a model-based algorithm in \cite{shi2024distributionally} is proposed by penalizing the robust value estimates with a carefully designed data-driven penalty term.
	
	The structure of this paper is organized as follows: Section 2 formulates the exploratory state dynamics equation under Knightian uncertainty, introduces the corresponding entropy-regularized relaxed stochastic control problem, and defines the value function. Section 3 employs stochastic analysis theory and the dynamic programming principle under nonlinear expectation to solve the HJB equation and derive the optimal randomized feedback control. Section 4 examines the LQ case, obtaining optimal randomized feedback control. Section 5 establishes the solvability equivalence between non-exploratory and exploratory LQ problems under Knightian uncertainty and analyzes exploration cost. Section 6 concludes the paper with a comprehensive summary.
	
	\section{Basic Model}
	\subsection{Preliminaries on Sublinear Expectations}
	%% Inline mathematics is tagged between $ symbols.
	In this paper, time is continuous and the uncertainty is characterized by a sublinear expectation space with filtration $(\Omega, \mathcal{H}, \{\mathcal{F}_t\}_{t \geq 0}, \hat{\mathbb{E}})$. Here, $\Omega$ represents the given set of elementary events, and a vector lattice on it $\mathcal{H}$ is the linear space consisting of real-valued functions defined on $\Omega$. $\mathcal{F}_t := \mathcal{B}(\Omega_t)$ represents the Borel \(\sigma\)-algebra over $\Omega_t$ and $\Omega_t := \left\{ \omega_{.\wedge t} : \omega \in \Omega \right\}$ denotes the set of all paths stopped at time $t$.
	$\hat{\mathbb{E}}$ denotes the sublinear expectation. Peng \cite{peng2010nonlinear} provided the following definition of sublinear expectation:
	\newtheorem{definition}{Definition}
	\begin{definition}
		A sublinear expectation $\hat{\mathbb{E}}$ is a (nonlinear) functional $\hat{\mathbb{E}}$: $\mathcal{H} \mapsto \mathbb{R}$ defined on a space $\mathcal{H}$ of random variables satisfying the following properties: 
		\begin{enumerate}
			\item[(i)] \textbf{Monotonicity}: For all random variables $X, Y \in \mathcal{H}$ that satisfy $X \geq Y$, we have $\hat{\mathbb{E}}[X] \geq \hat{\mathbb{E}}[Y]$;
			\item[(ii)] \textbf{Constant preserving}: $\hat{\mathbb{E}}[c] = c$, for $c \in \mathbb{R}$;
			\item[(iii)] \textbf{Sub-additivity}: For $\forall X, Y \in \mathcal{H}$, $\hat{\mathbb{E}}[X + Y] \leq \hat{\mathbb{E}}[X] + \hat{\mathbb{E}}[Y]$;
			\item[(iv)] \textbf{Positive homogeneity}: $\hat{\mathbb{E}}[aX] = a\hat{\mathbb{E}}[X]$, for  $\forall X\in \mathcal{H}$, $a \geq 0$.
		\end{enumerate}
	\end{definition}
	Let $\mathcal{E}$ denote the lower expectation, i.e., for any random variable $X$, there is $\mathcal{E}[X] = -\hat{\mathbb{E}}[-X]$.
	
	Under the given sublinear expectation framework, we consider RL-based optimal control problems. The RL agent interacts with and learns from an unknown environment through exploration, directly outputting an optimal or near-optimal control policy. In this process, the agent operates in complex and dynamic environments that may be subject to uncertain disturbances. These disturbances are multifaceted and vary across different real-world scenarios. For example, in financial markets, when an agent explores and learns to output an optimal investment policy, uncertainties arise from factors such as a company’s internal operations, competition among firms, and macroeconomic policies. Since the probability distribution of such uncertainties cannot be precisely determined by a probability distribution, they are referred to as Knightian uncertainty. The dynamics of the controlled system are often driven by $G$-Brownian motion $\{ B(t), t \geq 0 \}$ under $G$-expectation. The so-called $G$-Brownian motion is likewise clearly defined by Peng \cite{peng2010nonlinear}. Notably, the stability of systems driven by $G$-Brownian motion—particularly in the presence of time delays—has been rigorously analyzed; for instance, Fei et al. \cite{fei2025quasi} established quasi-sure exponential stability criteria for stochastic differential delay systems driven by G-Brownian motion (SDDE-GBM).
	
	And $B(1)$ obeys the $G$-normal distribution $\mathcal{N}(0, [\underline{\sigma}^2, \overline{\sigma}^2])$, where $\underline{\sigma}$ and $\overline{\sigma}$ are the minimum and maximum volatilities under Knightian uncertainty, respectively.
	
	Although the aforementioned $G$-normal distribution $\mathcal{N}(0, [\underline{\sigma}^2, \overline{\sigma}^2])$ is theoretically abstract, its characterizing parameters---the lower variance $\underline{\sigma}^2$ and upper variance $\overline{\sigma}^2$---can be estimated from empirical data using the $\varphi$-max-mean algorithm proposed by Peng \cite{peng2017}. Specifically, the i.i.d. sample data are partitioned into $m$ batches. For each batch $k$, a sample variance $\sigma_k^2$ can be computed. The optimal asymptotically unbiased estimators for the lower and upper variances are then given by
	\[
	\hat{\underline{\sigma}}^2 := \min_{1 \leq k \leq m} \sigma_k^2, \quad \hat{\overline{\sigma}}^2 := \max_{1 \leq k \leq m} \sigma_k^2,
	\]
	where
	\[
	\sigma_k^2 := \frac{1}{n} \sum_{j=1}^{n} (x_{n(k-1)+j} - \mu_k)^2,
	\]
	$x_{n(k-1)+j}$ represents the $j$-th data point in the $k$-th batch, and $\mu_k$ signifies the sample mean of the $k$-th batch.
	
	Denote $
	G(\alpha) := \frac{1}{2} \hat{\mathbb{E}}[\alpha B(1)^2] = \frac{1}{2} (\overline{\sigma}^2 \alpha^+ - \underline{\sigma}^2 \alpha^-),
	$
	where $\alpha^{+} = \max(\alpha, 0)$ and $\alpha^{-} = \max(-\alpha, 0)$ denote the positive and negative parts of $\alpha$, respectively. Here,
	\[
	\hat{\mathbb{E}}[B(1)^2] = \overline{\sigma}^2, \quad \mathcal{E}[B(1)^2] = \underline{\sigma}^2, \quad 0 < \underline{\sigma} \leq \overline{\sigma} < \infty.
	\]
	
	A property holds ``quasi-surely'' (q.s.) if there exists a polar set $D$ with $\mathbb{V}(D) = 0$ such that the property holds for any $\omega \in \Omega \setminus D$. Where upper capacity $\mathbb{V}$
	are defined as: 
	\[
	\mathbb{V}(A) := \sup_{P \in \mathcal{P}} P(A), \quad \forall A \in \mathcal{F}.
	\]

	The goal of RL is for the agent to learn an optimal policy that maximizes cumulative rewards. In practical applications such as robotic control and autonomous driving, the agent must make reliable decisions in complex, dynamic environments. This necessitates robust policies, which motivates our assumption that the RL agent is ambiguity aversion.
	
	\subsection{RL Modeling in Continuous Time and Space under Knightian Uncertainty}
	
	Assume that the $G$-Brownian motion $B = \{ B(t),\ t \geq 0 \}$ is a scalar process. Given an ``action space'' $U$, which represents constraints on the decision (``control'' or ``action'') of the agent. An admissible (open-loop) control $u = \{u_t, t \geq 0\}$ is an $\{\mathcal{F}_t\}_{t \geq 0}$-adapted measurable process taking values in $U$.
	
	The current state of the system satisfies the following stochastic differential equation:
	\begin{equation}\label{eq:1}
		dx_t^u = b \left( x_t^u, u_t \right) dt + \sigma \left( x_t^u, u_t \right) dB_t, t > 0; x_0^u = x \in \mathbb{R}.
		\tag{1}
	\end{equation}
	In this paper, $x$ is the generic variable that represents the current state of the system. The control objective is to maximize the expected total discounted return expressed by the following value function:
	\begin{equation}\label{eq:2}
		V^{\text{ne}}(x) := \sup_{u \in \mathcal{A}^{\text{ne}}(x)} \mathcal{E} \left[ \int_0^\infty e^{-\rho t} r(x_t^u, u_t) dt \bigg| x_0^u = x \right],
		\tag{2}
	\end{equation}
	where $r$ is the reward function, $\rho > 0$ is the discount rate, $\mathcal{A}^{\text{ne}}(x)$ denotes the set of all admissible controls which in general may depend on $x$, and $V^{\text{ne}}(x)$ denotes the non-exploratory objective value function.
	
	\begin{remark}
		In this definition, the lower expectation operator $\mathcal{E}[\cdot]$ is employed based on the assumption that the RL agent is ambiguity-averse. When confronted with an environment of unknown probabilistic model, an ambiguity-averse agent naturally adopts a robust decision-making criterion, i.e., by considering the worst-case scenario. The lower expectation $\mathcal{E}[\cdot]$, which computes the infimum of expected values over a set of plausible probability measures, mathematically captures this conservative perspective and provides a robust performance guarantee. The value function $V^{\text{ne}}(x)$ represents the cumulative reward. For an ambiguity-averse agent, using $\mathcal{E}[\cdot]$ to compute this expected value provides a conservative estimate, i.e., the cumulative reward that the agent can secure even under the most adverse probability model, thereby establishing a robust benchmark for policy evaluation.
		
		In contrast, using the upper expectation $\hat{\mathbb{E}}[\cdot]$ implies an ambiguity-loving agent who is optimistic and believes the environment will evolve in the most favorable way. This contradicts our core hypothesis of an ambiguity-averse agent seeking robustness. Consequently, the use of $\hat{\mathbb{E}}[\cdot]$ would lead to fundamentally different results and interpretations: the value function would represent an optimistic assessment of cumulative rewards, and the derived optimal policy would lack the desired robustness properties against model uncertainty.
	\end{remark}

	In the model-known case, the control $u_t$ in equation \eqref{eq:1} is deterministic, whereas in the RL context, the underlying model is unknown, meaning that the functions $b$, $\sigma$ are unspecified and thus requires dynamic learning. The agent interacts with and learns from an unknown environment through exploration. The key idea involves modeling exploration via randomized control $\theta = \{\theta_t(u), t \geq 0\}$ over the control space $U$, where each ``trial'' is sampled from $U$. Consequently, we can extend the classical control concept to incorporate randomized control. The agent executes a control for $N$ rounds over the same time horizon. In each round, a classical control is sampled from the randomized control $\theta$. The reward of such a policy becomes sufficiently precise when $N$ is large. Therefore, in order to evaluate this policy distribution in a continuous time setting, the limiting case of $N \to \infty$ needs to be considered.
	
	For the general reward function $r \left( x_t^u, u_t \right)$ that depends on both state and action, we first need to characterize how exploration modifies the state dynamics equation \eqref{eq:1} through an appropriately ``exploratory'' version. To achieve this, we investigate the effects of $N$ rounds of repeated learning under a given randomized control policy (e.g., $\theta$). Let $B_t^i$, $i = 1, 2, \dots, N$ be $N$ independent sample paths of the $G$-Brownian motion, and let $x_t^i$, $i = 1, 2, \dots, N$ denote the corresponding copies of the state process under controls $u^i$, $i = 1, 2, \dots, N$, respectively, where each $u^i$ is sampled from $\theta$. Then, for $i = 1, 2, \dots, N$, the increments of these state processes satisfy
	\begin{equation}\label{eq:3}
		\scalebox{0.9}{%
			$\displaystyle
			\Delta x_t^i \equiv x_{t+\Delta t}^i - x_t^i
			\approx b(x_t^i, u_t^i)\Delta t + \sigma(x_t^i, u_t^i)\bigl(B_{t+\Delta t}^i - B_t^i\bigr).
			$}
		\tag{3}
	\end{equation}
	
	Each such process $x^i$, $i = 1, 2, \dots, N$ can be viewed as an independent sample drawn from the exploratory state dynamics $X^\theta$. The superscript $\theta$ in $X^\theta$ indicates that each sample path $x^i$ is generated according to the dynamics \eqref{eq:3}, with the corresponding control $u^i$ being independently sampled under the randomized policy $\theta$.
	
	Based on dynamics \eqref{eq:3}, when \( N \to \infty \), we deduce
	\[
	\begin{aligned}
		\frac{1}{N} \sum_{i=1}^N \Delta x_t^i 
		&\approx \frac{1}{N} \sum_{i=1}^N b(x_t^i, u_t^i) \Delta t \\
		&\quad + \frac{1}{N} \sum_{i=1}^N \sigma(x_t^i, u_t^i) \left( B_{t+\Delta t}^i - B_t^i \right).
	\end{aligned}
	\]
	Let $\bar{b} = \frac{1}{N} \sum\limits_{i=1}^N b(x_t^i,u_t^i) $, with  $\lim\limits_{N \to \infty} \bar{b} \in [\underline{\mu}^b, \bar{\mu}^b]$. Here,
	\[
	\underline{\mu}^b = \mathcal{E} \left[ \int_U b(X_t^\theta, u) \theta_t(u) du \right],
	\]
	\[
	\overline{\mu}^b = \hat{\mathbb{E}} \left[ \int_U b(X_t^\theta, u) \theta_t(u) du  \right].
	\]
	According to the strong law of large numbers for sublinear expectations in (Chen \cite{Chen2016}, Theorem 3.1), $\mathcal{E}(B_{t+\Delta t}^i - B_t^i) = 0$ and $\hat{\mathbb{E}}(B_{t+\Delta t}^i - B_t^i) = 0$, we derive
	\begin{equation}\label{eq:4}
		\begin{split}
			\frac{1}{N} \sum_{i=1}^N \Delta x_t^i \xrightarrow{\text{q.s.}}
			\bigg[ & \mathcal{E} \left[ \int_U b(X_t^\theta, u) \theta_t(u) du \Delta t \right], \\
			& \hat{\mathbb{E}} \left[ \int_U b(X_t^\theta, u) \theta_t(u) du \Delta t \right] \bigg].
		\end{split}
		\tag{4}
	\end{equation}
	
	In the above, we have implicitly applied the (reasonable) assumption that both \(\theta_t\) and \(X_t^\theta\) are independent of the increments \(G\)-Brownian motion sample paths, and that these increments are identically distributed over the time interval \([t,t+\Delta t]\). Similarly, when \(N \to \infty\), we have
	\begin{align*}
		\frac{1}{N} \sum_{i=1}^{N} \left( \Delta x_t^i \right)^2 \approx \frac{1}{N} \sum_{i=1}^{N} \sigma^2(x_t^i, u_t^i) \left( B_{t+\Delta t}^i - B_t^i \right)^2.
	\end{align*}
	Moreover, we get
	\begin{equation}\label{eq:5}
		\scalebox{0.9}{%
			$\begin{aligned}
				\frac{1}{N} \sum_{i=1}^{N} (\Delta x_t^i)^2 \xrightarrow{\text{q.s.}}
				\bigg[ & \mathcal{E}\!\left[ \int_U \sigma^2(X_t^\theta,u)\theta_t(u)\,du\,d\langle B\rangle_t \right],\\
				& \hat{\mathbb{E}}\!\left[ \int_U \sigma^2(X_t^\theta,u)\theta_t(u)\,du\,d\langle B\rangle_t \right]
				\bigg].
			\end{aligned}$%
		}
		\tag{5}
	\end{equation}
	
	As we have seen, both \(\Delta x_t^i\) and \(\left(\Delta x_t^i\right)^2\) are affected by repeated learning given the policy \(\theta\).
	
	Formulas \eqref{eq:4} and \eqref{eq:5} motivates us to propose the exploratory version of state dynamics under Knightian uncertainty, namely,
	\begin{equation}\label{eq:6}
		dX_{t}^{\theta} = \tilde{b}(X_{t}^{\theta}, \theta_{t})dt + \tilde{\sigma}(X_{t}^{\theta}, \theta_{t})dB_{t}, t > 0; X_{0}^{\theta} = x \in \mathbb{R}
		\tag{6}.
	\end{equation}
	
	The coefficients \( \tilde{b}(\cdot, \cdot) \) and $\tilde{\sigma}(\cdot, \cdot)$ are the exploratory drift term and the exploratory volatility term, respectively, and they are defined as follows:
	\begin{equation}\label{eq:7}
		\tilde{b}(y, \theta) := \int_{U} b(y, u) \theta(u) du, \quad y \in \mathbb{R}, \theta \in \mathcal{P}(U)
		\tag{7},
	\end{equation}
	\begin{equation}\label{eq:8}
		\tilde{\sigma}(y, \theta) := \sqrt{\int_{U} \sigma^{2}(y, u) \theta(u) du}, \quad y \in \mathbb{R}, \theta \in \mathcal{P}(U)
		\tag{8},
	\end{equation}
	where \( \mathcal{P}(U) \) is the set of density functions of probability measures on \( U \) that are absolutely continuous with respect to the Lebesgue measure.
	
	Similarly, according to the strong law of large numbers for sublinear expectations (Chen \cite{Chen2016}, Theorem 3.1), when \( N \to \infty \), we deduce
	\[
	\scalebox{0.9}{$
		\begin{aligned}
			\frac{1}{N} \sum_{i=1}^N e^{-\rho t} r(x_t^i, u_t^i) \Delta t \xrightarrow{\mathrm{q.s.}}
			\bigg[ & \mathcal{E}\!\left[ e^{-\rho t}\! \int_U r(X_t^\theta, u)\theta_t(u)\,du\,\Delta t \right],\\
			& \hat{\mathbb{E}}\!\left[ e^{-\rho t}\! \int_U r(X_t^\theta, u)\theta_t(u)\,du\,\Delta t \right]
			\bigg].
		\end{aligned}
		$}
	\]
	Therefore, \( r \) in the \eqref{eq:2} needs to be modified to be the exploratory reward
	\[
	\tilde{r}(y, \theta) := \int_U r(y, u) \theta(u) du, \quad y \in \mathbb{R}, \theta \in \mathcal{P}(U).
	\]
	
	When there is no Knightian uncertainty (i.e., $\underline{\sigma} = \overline{\sigma} = \sigma$), $d\langle B \rangle_t = \sigma^2 dt$, and the entire derivation reduces to the classical framework of Wang et al. \cite{wang2020reinforcement}.

	\subsection{Entropy Regularization}
	We introduce a relaxed stochastic control formulation to model the exploration and exploitation in RL. In fact, if the model were fully known, there would be no need for exploration and learning at all, and the randomized controls would all degenerate to Dirac measures—this would place us strictly within the domain of stochastic control under nonlinear expectation. Therefore, in the RL context, it is necessary to incorporate a “regularization term” to account for model parameter uncertainty and encourage exploration. Here, we use Shannon’s differential entropy to quantify the level of exploration
	\[
	\mathcal{H}(\theta) := -\int_{U} \theta(u) \ln \theta(u) du, \quad \theta \in \mathcal{P}(U).
	\]
	
	Introduce the following entropy-regularized relaxed stochastic control problem in the Knightian uncertainty setting
	\begin{equation}\label{eq:9}
		\begin{split}
			V(x) &:= \sup_{\theta \in \mathcal{A}(x)} \mathcal{E} \bigg[ \int_0^\infty e^{-\rho t} \bigg( \int_U r \left( X_t^\theta, u \right) \theta_t (u) du \\
			&\quad -\lambda \int_U \theta_t (u) \ln \theta_t (u) du \bigg) dt \bigg| X_0^\theta = x \bigg],
		\end{split}
		\tag{9}
	\end{equation}
	where \(\lambda > 0\) is the exogenous exploratory weight parameter capturing the trade-off between exploitation (the original reward function) and exploration (the entropy). \(\mathcal{A}(x) \) is the set of the admissible randomized controls (usually dependent on \(x\)), and \(V(x)\) denotes the exploratory objective value function.
	
	The precise definition of \(\mathcal{A}(x) \) depends on the specific dynamic model under consideration and the particular problem to be solved, which may vary from case to case. Here, we first outline some ``minimal'' requirements for \(\mathcal{A}(x) \). Let \( \mathcal{B}(U) \) denote the Borel \(\sigma\)-algebra on \( U \). An admissible randomized control is a measure-valued (or density-function-valued) process \( \theta = \{\theta_t, t \geq 0\} \) that satisfies at least the following properties:
	
	\begin{enumerate}
		\item[(i)] for each \( t \geq 0 \), \( \theta_t \in \mathcal{P}(U) \) q.s.;
		
		\item[(ii)] for each \( A \in \mathcal{B}(U) \), \(\left\{ \int_A \theta_t(u)du, t \geq 0 \right\}\) is \( \mathcal{F}_t \)-progressively measurable;
		
		\item[(iii)] if \( \theta \) is applied, then the stochastic differential equation \eqref{eq:6} will have a unique strong solution \( X^\theta = \{X_t^\theta, t \geq 0\} \);
		
		\item[(iv)] the lower expectation on the right-hand side of \eqref{eq:9} is finite.
	\end{enumerate}
	
	Finally, similar to the classical control case,  \(\mathcal{A}(x) \) contains open-loop randomized controls which are measure-valued stochastic processes. We will also consider randomized feedback controls. Specifically, a deterministic mapping \(\boldsymbol{\theta}(\cdot;\cdot)\) is called randomized feedback control subject to the following conditions:
	\begin{enumerate}
		\item[(i)] for each \( x \in \mathbb{R} \),  \(\boldsymbol{\theta}(\cdot;\cdot)\) is a density function;
		
		\item[(ii)] the dynamics of the system after applying the randomized feedback control \(\boldsymbol{\theta}(\cdot;\cdot)\) as follows has a unique strong solution \(\{X_t, t \geq 0\}\)
		\[
		dX_t = \tilde{b}(X_t, \boldsymbol{\theta}(\cdot; X_t))dt + \tilde{\sigma}(X_t, \boldsymbol{\theta}(\cdot; X_t))dB_t, t > 0,
		\]
		where \(X_0 = x \in \mathbb{R}\).
		\item [(iii)] the open-loop randomized control \(\theta = \{\theta_t, t \geq 0\} \in \mathcal{A}(x)\), where \(\theta_t = \boldsymbol{\theta}(\cdot; X_t)\).
	\end{enumerate}
	
	In this case, it can be said that the open-loop randomized control \( \theta \) is generated from the randomized feedback control \(\boldsymbol{\theta}(\cdot;\cdot)\) with respect to \( x \).
	\section{HJB Equations and Optimal Randomized Control}
	
	Below, we outline the general procedure for solving the optimization problem \eqref{eq:9}. We have
	\[
	\begin{split}
		V(x) &=\sup_{\theta \in \mathcal{A}(x)} \mathcal{E} \bigg[ \int_0^s e^{-\rho t} \left( \tilde{r} \left( X_t^\theta, \theta_t \right) + \lambda \mathcal{H} \left( \theta_t \right) \right) dt \\
		&\quad + e^{-\rho s} V \left( X_s^\theta \right) \bigg| X_0^\theta = x \bigg], s > 0.
	\end{split}
	\]
	\begin{theorem}
		Let $v$ be the general unknown solution to the following HJB equation:
		\begin{align*}
			\rho v(x) &= \max_{\theta \in\mathcal{P}(U)} \bigg( \tilde{r}(x, \theta) - \lambda \int_U \theta(u) \ln \theta(u) du \\
			&\quad +\tilde{\sigma}^2 (x, \theta) \widetilde{G}[v''(x)] + \tilde{b}(x, \theta) v'(x) \bigg), \quad x \in \mathbb{R},
		\end{align*}
		equivalently expressed as:
		\begin{equation}\label{eq:10}
			\begin{split}
				\rho v(x) &= \max_{\theta \in\mathcal{P}(U)} \int_U \bigg( r(x, u) - \lambda \ln \theta(u) \\
				&\quad +\sigma^2 (x, u) \widetilde{G}[v''(x)] + b(x, u) v'(x) \bigg) \theta(u) du,
			\end{split}
			\tag{10}
		\end{equation}
		where
		\[
		\widetilde{G}[v''(x)] := \frac{1}{2} \mathcal{E}[v''(x) B(1)^2] = \frac{1}{2} (\underline{\sigma}^2 [v''(x)]^+ - \overline{\sigma}^2 [v''(x)]^-).
		\]
		$ \theta \in \mathcal{P}(U) $, if and only if $
		\int_U \theta(u) du = 1, \, \theta(u) \geq 0 \text{ a.e. on } U.$
		
		The constrained optimization problem in \eqref{eq:10} admits a unique solution given by the randomized feedback control
		\begin{equation}\label{eq:11}
			\begin{aligned}
				\boldsymbol{\theta}^*(u;x) = \frac{ \exp\left( \frac{1}{\lambda} \Psi(x,u) \right) }{ \int_{U} \exp\left( \frac{1}{\lambda} \Psi(x,u) \right) \mathrm{d}u },
			\end{aligned}
			\tag{11}
		\end{equation}
		where $\Psi(x,u) = r(x,u) + \sigma^2(x,u) \widetilde{G}[v''(x)] + b(x,u)v'(x)$.
	\end{theorem}
	\begin{proof}
		By employing the nonlinear dynamical programming principle (Fei \cite{chen2021optimal}), we derive that
		the HJB equation \eqref{eq:10}.
		
		Treating the right-hand side of the HJB equation \eqref{eq:10} as a maximization problem for \( \theta(u) \), subject to the constraints \( \int_{U} \theta(u) d u=1 \) and \( \theta(u) \geq 0 \), we use the Lagrange multiplier method to construct the Lagrangian function
		\begin{equation*}
			\small
			\begin{split}
				L &= \int_{U} \big[ r(x, u) \theta(u) - \lambda \theta(u) \ln \theta(u) + \sigma^{2}(x, u) \theta(u) \widetilde{G}\left[v''(x)\right] \\
				&\quad+ b(x, u) \theta(u) v'(x) \big] du + \eta \left(1 - \int_{U} \theta(u) du\right),
			\end{split}
		\end{equation*}
		where \( \eta \) is the Lagrange multiplier.
		
		Take the derivative of \( L \) with respect to \( \theta(u) \) and set it to zero, we get
		\begin{equation*}
			\begin{split}
				\frac{\partial L}{\partial \theta(u)} &=  r(x, u) - \lambda (1 + \ln \theta(u)) + \sigma^{2}(x, u) \widetilde{G}\left[v''(x)\right] \\
				&\quad + b(x, u) v'(x) - \eta =0.
			\end{split}
		\end{equation*}
		Denote the constant term \( C = \exp\left(-1 - \frac{\eta}{\lambda}\right) \). Then, we have
		\begin{equation*}
			\begin{split}
				\theta(u) &= C \cdot \exp\bigg(\frac{1}{\lambda}\big[  r(x, u) + \sigma^{2}(x, u) \widetilde{G}\left[v''(x)\right] \\
				&\quad + b(x, u) v'(x)\big]\bigg).
			\end{split}
		\end{equation*}
		From the normalization condition of the probability density \( \int_{U} \theta(u) d u = 1 \), we obtain
		\begin{equation*}
			\begin{split}
				C &= \bigg[ \int_U \exp\bigg( \frac{1}{\lambda} \Big[  r(x,u) + \sigma^2(x,u)\widetilde{G}[v''(x)] \\
				&\quad + b(x,u)v'(x) \Big] \bigg) du \bigg]^{-1}.
			\end{split}
		\end{equation*}
		Substituting $C$ into the expression for \( \theta(u) \), we obtain equation \eqref{eq:11}. The proof hence is complete.
	\end{proof}
	
	The optimal randomized feedback control $\boldsymbol{\theta}^*(u;x)$ in \eqref{eq:11} exhibits case-dependent behavior through the operator $\widetilde{G}[v''(x)]$, a direct manifestation of the Knightian uncertainty. This stands in sharp contrast to the characteristic obtained by Wang et al. \cite{wang2020reinforcement} under the classical linear expectation scenario, where the variance depends solely on a fixed parameter $\sigma^2$. Specifically:
	
	When $v''(x) > 0$, the control depends on the lower variance $\underline{\sigma}^2$. This reflects the agent's prudent attitude in a potential ``gain" region (where the value function is convex): even if favorable uncertainty exists, it adopts the most conservative volatility estimate to avoid risk.
	
	When $v''(x) < 0$, the control depends on the upper variance $\overline{\sigma}^2$. This profoundly illustrates the agent's extreme pessimism in a potential ``loss" region (where the value function is concave): due to ambiguity aversion, the agent tends to amplify risks, employing the worst-case volatility estimate to assess downside risk, thereby leading to stronger risk aversion and suppression of exploration.
	
	When $v''(x) = 0$ (and $v'(x) \neq 0$), it reduces to containing only the drift term.
	
	If the agent is ambiguity-loving, then the behavior of the optimal randomized feedback control $\boldsymbol{\theta}^*(u;x)$ in \eqref{eq:11} exhibits case-dependence through the operator ${G}[v''(x)]$ that is opposite to the ambiguity-averse case described above.
	
	When $v''(x) > 0$, the control depends on the upper variance $\overline{\sigma}^2$, reflecting the agent's optimistic attitude in potential gain regions and amplifying exploration to pursue the best possible outcomes.
	
	When $v''(x) < 0$, the control depends on the lower variance $\underline{\sigma}^2$, demonstrating the agent's tendency to underestimate downside risks in potential loss regions.
	
	The above \eqref{eq:11} contributes to the qualitative understanding of optimal exploration in Knightian uncertainty environment. We will investigate this further in the next section.
	
	\section{The Linear-Quadratic Case}
	
	We now focus on a class of entropy-regularized (relaxed) stochastic control problems with linear state dynamics and quadratic rewards under Knightian uncertainty, in which
	\begin{equation}\label{eq:12}
		b(x, u) = Ax + Fu, \quad \sigma(x, u) = Cx + Du, \quad x, u \in \mathbb{R},
		\tag{12}
	\end{equation}
	where \( A, F, C, D \in \mathbb{R} \). The reward function is given by
	\begin{equation}\label{eq:13}
		\scalebox{0.9}{$
			r(x, u) = -\left( \frac{M}{2}x^2 + Ixu + \frac{K}{2}u^2 + Px + Qu \right), \quad x, u \in \mathbb{R},
			$}
		\tag{13}
	\end{equation}
	where \( M \geq 0, K > 0, I, P, Q \in \mathbb{R} \). As is standard in LQ control theory, we assume the control set is unconstrained, i.e., \( U = \mathbb{R} \).
	
	Fix the initial state \( x \in \mathbb{R} \). The mean and variance processes \( \mu_t, \sigma_t^2, t \geq 0 \) of any open-loop randomized control \( \theta \in \mathcal{A}(x) \) are characterized as follows:
	\[
	\mu_t := \int_{\mathbb{R}} u\theta_t(u)du, \quad \sigma_t^2 := \int_{\mathbb{R}} u^2\theta_t(u)du - \mu_t^2.
	\]
	
	The stochastic differential equation \eqref{eq:6} can be expressed as:
	\begin{equation}\label{eq:14}
		\scalebox{0.87}{%
			$\begin{aligned}
				dX_t^\theta
				&= \int_U b(X_t^\theta, u)\theta_t(u)\,du\,dt
				+ \sqrt{\int_U \sigma^2(X_t^\theta, u)\theta_t(u)\,du}\,dB_t \\[2pt]
				&= (AX_t^\theta + F\mu_t)\,dt
				+ \sqrt{(C X_t^\theta + D\mu_t)^2 + D^2\sigma_t^2}\,dB_t,
			\end{aligned}$%
		}
		\tag{14}
	\end{equation}
	where \(t > 0, X_0^\theta = x \nonumber\).
	
	Further define
	\[
	g(X_t^\theta, \theta_t) := \int_{\mathbb{R}} r(X_t^\theta, u)\theta_t(u)du - \lambda \int_{\mathbb{R}} \theta_t(u)\ln\theta_t(u)du.
	\]
	Next, we define the associated admissible control set $\mathcal{A}(x)$: $\theta \in \mathcal{A}(x)$ if
	
	\begin{enumerate}
		\item[(i)] for each $t \geq 0$, $\theta_t \in \mathcal{P}(U)$ q.s.;
		
		\item[(ii)] for each $A \in \mathcal{B}(U)$, $\left\{ \int_A \theta_t(u)du, t \geq 0 \right\}$ is $\mathcal{F}_t$-progressively measurable;
		
		\item[(iii)]for each $t \geq 0$, $\hat{\mathbb{E}}\left[\int_0^t (\mu_s^2 + \sigma_s^2) ds\right] < \infty$;
		
		\item[(iv)]for $\left\{ X_t^\theta, t \geq 0 \right\}$ solving \eqref{eq:14}, \\$\liminf_{T \to \infty} e^{-\rho T} \hat{\mathbb{E}}\left[(X_T^\theta)^2\right] = 0$;
		
		\item[(v)] for $\left\{ X_t^\theta, t \geq 0 \right\}$ solving \eqref{eq:14}, \\$\hat{\mathbb{E}}\left[\int_0^\infty e^{-\rho t} |g(X_t^\theta, \theta_t)| dt\right] < \infty$.
	\end{enumerate}
	
	In the above, condition (iii) is to ensure that for any \(\theta \in \mathcal{A}(x)\), the drift and volatility terms in \eqref{eq:14} satisfy the global Lipschitz condition and a class of linear growth conditions in the state variables, and hence, based on Peng \cite{peng2010nonlinear}, the SDE \eqref{eq:14} has a unique strong solution \(X^\theta\). Condition (iv) is to ensure that the nonliner dynamic programming principle as well as the verification theorem are applicable for this model. Finally, condition (v) ensures that the reward is finite.
	
	We now introduce the entropy-regularized relaxed LQ control problem under Knightian uncertainty
	\begin{equation}\label{eq:15}
		\begin{split}
			V(x) &= \sup_{\theta \in \mathcal{A}(x)} \mathcal{E}\bigg[\int_0^\infty e^{-\rho t}\Big(\int_{\mathbb{R}} r(X_t^\theta, u)\theta_t(u)du \\
			&\quad-\lambda\int_{\mathbb{R}}\theta_t(u)\ln\theta_t(u)du\Big)dt\,\bigg|\,X_0^\theta = x\bigg].
		\end{split}
		\tag{15}
	\end{equation}
	
	Based on the derivation in the previous section, substituting equations \eqref{eq:12} and \eqref{eq:13} into \eqref{eq:11}, through straightforward calculation under the condition that $K - 2D^2\widetilde{G}[v''(x)] > 0$, we obtain the optimal randomized feedback control in \eqref{eq:11} reduces to
	\begin{equation}\label{eq:16}
		\tag{16}
		\boldsymbol{\theta}^*(u; x) = \frac{\exp \left( - \left( u -  \mu_1 \right)^2 \big/ {2\sigma_1^2} \right)}{\int_{\mathbb{R}} \exp \left( - \left( u - \mu_1  \right)^2 \big/ {2\sigma_1^2} \right) du},
	\end{equation}
	where
	\[
	\mu_1 = \frac{2CDx\widetilde{G}[v''(x)] + Fv'(x) -Ix- Q}{K - 2D^2\widetilde{G}[v''(x)]},
	\]
	and
	\[
	\sigma_1^2 = \frac{\lambda}{K - 2D^2\widetilde{G}[v''(x)]}.
	\] 
	In this case, the optimal randomized feedback control follows a Gaussian distribution.
	
	The mean \( \mu_1 \) represents the central tendency of the control policy \( \theta \), indicating the optimal action the agent tends to choose in the current state. This reflects the agent's optimal utilization of known information in the current state. We observe that for the general case where the reward depends on both state and control, the optimal action $(\mu_1)$ on $\mathbb{R}$ also depends on the current state $x$.
	
	The variance \( \sigma_1^2 \) characterizes the dispersion of the control policy \( \theta \), representing the randomness in the agent's exploration of unknown actions. The larger the variance, the wider the scope of the agent's exploration, while the smaller the variance, the more the agent's preference for actions closer to the mean, i.e., the current optimal action.
	
	Notice that the condition \( K - 2D^2 \widetilde{G}[v''(x)] > 0,\ x \in \mathbb{R} \) is fundamental to ensuring the well-posedness of the Gaussian policy in \eqref{eq:16}. This condition serves not merely as a theoretical assumption but as a verification criterion that must be assessed based on the estimated bounds of Knightian uncertainty, $\hat{\underline\sigma}^2$ and $\hat{\overline\sigma}^2$, which can be obtained from empirical data using the $\varphi$-max-mean algorithm (Peng \cite{peng2017} and Deng et al. \cite{deng2019}).
	
	When $v''(x) = 0$, the condition simplifies to $ K > 0 $, which clearly holds. The optimal randomized feedback control is Gaussian distribution.
	
	When $v''(x) > 0$, we have $\widetilde{G}[v''(x)] = \frac{1}{2} \underline{\sigma}^2 v''(x)$. The condition becomes $
	K - D^2 \underline{\sigma}^2 v''(x) > 0.$
	Given the estimate $\hat{\underline\sigma}^2$, the agent can conservatively verify a sufficient condition for a well-defined Gaussian policy: $
	K - D^2 \hat{\underline\sigma}^2v''(x) > 0.$
	
	The analysis reveals that an environment with a higher estimated $\hat{\underline\sigma}^2$ (indicating a higher minimum level of uncertainty and thus a narrower Knightian uncertainty interval) leads to an increase in the exploration variance $\sigma_1^2$. This phenomenon occurs because
	the RL agent is ambiguity aversion: when Knightian uncertainty
	diminishes, the agent may perceive the environment as
	more controllable and consequently becomes more willing to
	try new actions (exploration) to further optimize its policy.
	
	When $v''(x) < 0$, we have $\widetilde{G}[v''(x)] = \frac{1}{2} \overline{\sigma}^2 v''(x)$. The condition $ K - D^2 \overline{\sigma}^2 v''(x) > 0 $ holds automatically given $ v''(x) < 0 $. An environment with a higher estimated $\hat{\overline\sigma}^2$ (indicating a higher degree of Knightian uncertainty) leads to a decrease in the exploration variance $\sigma_1^2$. This forces the agent to adopt a more conservative policy, reducing exploration and favoring the current optimal action to hedge against the elevated model uncertainty. The agent’s aversion to Knightian uncertainty implies it will avoid making decisions in highly uncertain environments.
	
	Our condition $K - 2D^2 \widetilde{G}[v''(x)] > 0$ is a generalization of the condition $N - D^2 v''(x) > 0$ in Wang et al. \cite{wang2020reinforcement} under Knightian uncertainty. The condition in \cite{wang2020reinforcement} arises in a stochastic environment with a known volatility $\sigma$, leading to a linear dependence on $v''(x)$. In contrast, our condition accounts for Knightian uncertainty in the volatility and introduces nonlinearity and asymmetry through the $G$-expectation.
	
	Substituting \eqref{eq:16} into \eqref{eq:10},  the HJB equation becomes 
	\begin{equation}\label{eq:17}
		\begin{split}
			\rho v(x) &= \frac{(2CDx\widetilde{G}[v''(x)] + Fv'(x) - Ix - Q)^2}{K - 2D^2\widetilde{G}[v''(x)]} 
			\\
			&\quad+ \frac{1}{2}(2C^2\widetilde{G}[v''(x)] - M)x^2 
			+ (Av'(x) - P)x\\
			&\quad+ \frac{\lambda}{2} \left[ \ln \left( \frac{2\pi e\lambda}{K - 2D^2\widetilde{G}[v''(x)]} \right) - 1 \right].
		\end{split}
		\tag{17}
	\end{equation}
	
	A smooth solution of the HJB \eqref{eq:17} is given by the following equation:
	\[
	v(x) = \frac{1}{2} k_2 x^2 + k_1 x + k_0,
	\]
	where
	
	\begin{center}
		$\begin{aligned}
			\rho k_2 &= \frac{2(2CD\widetilde{G}[v''(x)] + Fk_2 - I)^2}{K - 2D^2\widetilde{G}[v''(x)]} + 2C^2\widetilde{G}[v''(x)] + 2Ak_2 - M, \\
			\rho k_1 &= \frac{2(2CD\widetilde{G}[v''(x)] + Fk_2 - I)(Fk_1 - Q)}{K - 2D^2\widetilde{G}[v''(x)]} + Ak_1 - P, \\
			\rho k_0 &= \frac{(k_1F - Q)^2}{K - 2D^2\widetilde{G}[v''(x)]} + \frac{\lambda}{2} \left( \ln \left( \frac{2\pi e\lambda}{K - 2D^2\widetilde{G}[v''(x)]} \right) - 1 \right).
		\end{aligned}$
	\end{center}
	
	When $M$ is sufficiently large, it follows that $k_2 < 0$, and hence $v''(x) < 0$. In this case, we have
	\begin{equation}\label{eq:18}
		k_2 = \frac{2[(CD\overline{\sigma}^2 + F)k_2 - I]^2}{\rho(K - k_2D^2\overline{\sigma}^2)} + \rho[(C^2\overline{\sigma}^2 + 2A)k_2 - M], \tag{18}
	\end{equation}
	\begin{equation}\label{eq:19}
		k_1 = \frac{2[(CD\overline{\sigma}^2 + F)k_2 - I)(Fk_1 - Q)]}{\rho(K - k_2D^2\overline{\sigma}^2)} + \rho(Ak_1 - P), \tag{19}
	\end{equation}
	\begin{equation}\label{eq:20}
		\scalebox{0.94}{%
			$\displaystyle
			k_0 = \frac{(k_1F - Q)^2}{\rho(K - k_2D^2\overline{\sigma}^2)}
			+ \frac{\lambda}{2\rho}\left(
			\ln\!\left(\frac{2\pi e\lambda}{K - k_2D^2\overline{\sigma}^2}\right) - 1
			\right).
			$}
		\tag{20}
	\end{equation}
	
	Equation \eqref{eq:18} is a nonlinear equation in terms of $k_2$, which cannot be solved analytically but can be addressed using numerical methods (e.g., Newton's method, the bisection method, etc.). Once the numerical solution for $k_2$ is obtained, $k_1$ and $k_0$ can be computed analytically via direct substitution into formulas \eqref{eq:19} and \eqref{eq:20}.
	
	For this particular solution, it is given by $v(x)$ above. Since $M$ is large enough, $k_2 < 0$, $v$ is a concave function. This property is essential in proving that it is in fact the value function. On the other hand, $K - 2D^2\widetilde{G}[v''(x)]> 0$ (i.e., $K - k_2D^2\overline{\sigma}^2 > 0$) ensures that $k_0$ is well-defined.
	
	\begin{proposition}
		Assume that the reward function is given by \eqref{eq:13}.
		Then, the value function in \eqref{eq:15} is given by the following equation:
		\[
		V(x) = \frac{1}{2} k_2 x^2 + k_1 x + k_0, x \in \mathbb{R},
		\]
		where $k_2$, $k_1$ and $k_0$ are shown by the above \eqref{eq:18}, \eqref{eq:19}, \eqref{eq:20}, respectively. In addition, the optimal randomized feedback control is Gaussian distributed with a density function of
		\begin{equation}\label{eq:21}
			\boldsymbol{\theta}^*(u;x) = \mathcal{N} \left( u \big|\mu_2, \sigma_2^2 \right), \tag{21}
		\end{equation}
		where
		\[
		\mu_2 = \frac{(k_2(F + CD\overline{\sigma}^2) - I)x + k_1F - Q}{K - k_2D^2\overline{\sigma}^2},\]
		\[\sigma_2^2 = \frac{\lambda}{K - k_2D^2\overline{\sigma}^2}.
		\]
		Finally, the associated optimal state process $\{X_t^*, t \geq 0\}$ under $\boldsymbol{\theta}^*(\cdot; \cdot)$ is the unique solution
		of the following SDE:
		\begin{equation*}
			\begin{aligned}
				dX_t^* &= \bigg[ \Big( A + \tfrac{F(k_2(F + CD\overline{\sigma}^2) - I)}{K - k_2D^2\overline{\sigma}^2} \Big) X_t^* + \tfrac{F(k_1F - Q)}{K - k_2D^2\overline{\sigma}^2} \bigg] dt \\
				&\quad+ \bigg[ \Big( \big( C + \tfrac{D(k_2(F + CD\overline{\sigma}^2) - I)}{K - k_2D^2\overline{\sigma}^2} \big) X_t^* + \tfrac{D(k_1F - Q)}{K - k_2D^2\overline{\sigma}^2} \Big)^2 \\
				&\quad+ \tfrac{\lambda D^2}{K - k_2D^2\overline{\sigma}^2} \bigg]^{1/2} dB_t,
			\end{aligned}
		\end{equation*}
		where $X_0^*$ = x.
	\end{proposition}
	
	\begin{proof}
		Let $v(x) = \frac{1}{2} k_2 x^2 + k_1 x + k_0$ be the candidate value function with coefficients defined in \eqref{eq:18}-\eqref{eq:20}, which satisfies the HJB equation \eqref{eq:10}.
		
		Throughout this proof we fix the initial state \( x \in \mathbb{R} \). Let \( \theta \in \mathcal{A}(x) \) and \( X^\theta \) be the associated state process solving \eqref{eq:14} with \( \theta \) being used. Let \( T > 0 \) be arbitrary. Define the stopping times 
		\[
		\tau_n^\theta := \left\{ t \geq 0 : \int_0^t \left( e^{-\rho t} v'(X_t^\theta) \tilde{\sigma}(X_t^\theta, \theta_t) \right)^2 dt \geq n \right\},
		\] 
		for $n \geq 1$.
		
		Applying the $G$-Itô formula to $e^{-\rho t}v(X_t^\theta)$, we have
		\begin{align*}
			&e^{-\rho(T\wedge\tau_{n}^{\theta})}v(X_{T\wedge\tau_{n}^{\theta}}^{\theta}) \\
			&= v(x) + \int_{0}^{T\wedge\tau_{n}^{\theta}} e^{-\rho t} \bigg( -\rho v(X_{t}^{\theta})
			+ v'(X_t^\theta)\tilde{b}(X_t^\theta,\theta_t) \bigg) dt \\
			&\quad + \int_0^{T\wedge\tau_n^\theta} e^{-\rho t}v'(X_t^\theta)\tilde{\sigma}(X_t^\theta,\theta_t) dB_t\\
			&\quad+\frac{1}{2}\int_0^{T\wedge\tau_n^\theta}e^{-\rho t}v^{\prime\prime}(X_{t}^{\theta})\tilde{\sigma}^{2}(X_{t}^{\theta},\theta_{t})d\langle B \rangle_t.
		\end{align*}
		Using $d\langle B\rangle_t \leq \overline{\sigma}^2 dt$ and substituting the inequality derived from the HJB equation \eqref{eq:10}, we get
		
		\[
		\begin{aligned}
			-\rho v(X_{t}^{\theta}) 
			&\leq -\tilde{r}(X_{t}^{\theta}, \theta_t) + \lambda \mathcal{H}(\theta_t)  - \frac{1}{2}v^{\prime\prime}(X_{t}^{\theta})\tilde{\sigma}^{2}(X_{t}^{\theta},\theta_{t})\overline{\sigma}^2\\
			&\quad - \tilde{b}(X_{t}^{\theta}, \theta_t) v'(X_{t}^{\theta}).
		\end{aligned}
		\]
		Thus we have
		\begin{align*}
			&e^{-\rho(T\wedge\tau_{n}^{\theta})}v(X_{T\wedge\tau_{n}^{\theta}}^{\theta}) \\
			&\leq v(x) + \int_{0}^{T\wedge\tau_{n}^{\theta}} e^{-\rho t} \left(-\tilde{r}(X_t^\theta, \theta_t) \right. \\
			&\quad \left. + \lambda \int_{\mathbb{R}} \theta_t(u) \ln \theta_t(u) du\right) dt \\
			&\quad + \int_0^{T\wedge\tau_n^\theta} e^{-\rho t}v'(X_t^\theta)\tilde{\sigma}(X_t^\theta,\theta_t) dB_t.
		\end{align*}
		Taking expectation $\hat{\mathbb{E}}$ on both sides, we get
		\[
		\begin{aligned}
			&\hat{\mathbb{E}}\left[e^{-\rho(T \wedge \tau_n^\theta)} v(X_{T \wedge \tau_n^\theta}^\theta)\right] \\
			& \leq v(x) -\mathcal{E}\left[\int_{0}^{T \wedge \tau_n^\theta} e^{-\rho t} \left(\tilde{r}(X_t^\theta, \theta_t) \right. \right. \\
			&\quad \left. \left. - \lambda \int_{\mathbb{R}} \theta_t(u) \ln \theta_t(u) \, du \right) dt \right].
		\end{aligned}
		\]
		
		Since for any \( r \geq 2 \), there exists a constant \( C = C(r, T, \overline{\sigma}) \) such that
		$
		\hat{\mathbb{E}}\left[ \sup_{0 \leq t \leq T} |X(t)|^r \right] \leq C.
		$ (by Yin et al. \cite{yin2019quasi} and Deng et al. \cite{deng2025convergence})
		
		Letting $n \to \infty$, we deduce that
		\[
		\begin{aligned}
			&\hat{\mathbb{E}} \left[ e^{-\rho T} v(X_T^\theta) \right]\\
			& \leq v(x) - \mathcal{E}\left[\int_0^T e^{-\rho t} \left(\tilde{r}(X_t^\theta, \theta_t) \right.\right. \\
			&\quad \left.\left. - \lambda \int_{\mathbb{R}} \theta_t(u) \ln \theta_t(u) du\right) dt\right].
		\end{aligned}
		\]
		
		We recall the admissibility condition 
		$
		\liminf_{T \to \infty} e^{-\rho T} \mathbb{\hat{E}}\left[ (X_T^\theta)^2 \right] = 0.
		$
		Thus, together with the fact that \( k_2 < 0 \), leads to
		$
		\limsup_{T \to \infty} \mathbb{\hat{E}}\left[ e^{-\rho T} v(X_T^\theta) \right] = 0.
		$
		So, we have
		\[
		\begin{aligned}
			v(x) &\geq \mathcal{E}\left[\int_0^T e^{-\rho t} \left(\tilde{r}(X_t^\theta, \theta_t) 
			\right. \right.  \\
			& \left. \left. \quad- \lambda \int_{\mathbb{R}} \theta_t(u) \ln \theta_t(u)  du \right) dt \right],
		\end{aligned}
		\]
		for each \( x \in \mathbb{R} \) and \( \theta \in \mathcal{A}(x) \). Hence, \( v(x) \geq V(x) \), for all \( x \in \mathbb{R} \).
		
		Let $\theta^* = \{\theta_t^*, t \geq 0\}$ be the open-loop control distribution generated from the above feedback law along with the corresponding state process $\{X_t^*, t \geq 0\}$ with $X_0^* = x$, and
		assume for now that $\theta^* \in \mathcal{A}(x)$. 
		
		Since the candidate value function \( v(x) \) is a solution satisfying the HJB equation \eqref{eq:10} and \(\theta^*\) is the optimal policy, it is obvious that
		\[
		\begin{aligned}
			v(x) &\leq\mathcal{E}\left[\int_0^\infty e^{-\rho t} \left(\tilde{r}(X_t^*, \theta_t^*) \right.\right. \\
			&\quad \left.\left. - \lambda \int_{\mathbb{R}} \theta_t^*(u) \ln \theta_t^*(u) du\right) dt\right]\\
			&=V(x),
		\end{aligned}
		\]
		for any $x \in \mathbb{R}$. This proves that $v$ is indeed the value function, namely $v = V$.  
		
		That $\theta^* \in \mathcal{A}(x)$ follows from Lemma 1 in Proposition 2, which provides the necessary integrability and asymptotic conditions. This completes the proof of Proposition 1.
	\end{proof}
	
	Different from the proof by Wang et al. \cite{wang2020reinforcement} under the linear expectation framework, our proof is established within the sublinear expectation framework and requires the \textit{G}-Itô's formula as well as the subadditivity of sublinear expectation.
	
	\begin{remark}
		The discount rate $\rho$ plays a pivotal role in regulating the exploration variance of the optimal randomized policy, albeit indirectly. A higher discount rate $\rho$, indicating a more myopic agent, reduces the perceived value of future rewards. This diminishes the agent's incentive to explore for long-term benefits, resulting in a decrease in the optimal exploration variance $\sigma_2^2$. Conversely, a lower discount rate $\rho$, representing a far-sighted agent, increases the valuation of future outcomes. This motivates the agent to engage in more extensive exploration, thereby increasing $\sigma_2^2$. This modulation effect is channeled through the coefficient $k_2$ in the solution to the HJB equation \eqref{eq:17}, which is intrinsically dependent on $\rho$ (see \eqref{eq:18}). Thus, the discount rate $\rho$ serves as a crucial tuning parameter for the exploration-exploitation trade-off, complementing the direct effect of the exploratory weight $\lambda$.
	\end{remark}
	
	\section{Exploration Cost and Its Impact}
	
	Motivated by the necessity of exploration in RL environments with typical unknown dynamics, we formulate and analyze a novel class of stochastic control problems under Knightian uncertainty that incorporates both entropy-regularized and relaxed control. We further derive closed-form solutions and provide verification theorems for an important class of LQ problems. In this section, we will conduct a comparative analysis between exploratory and non-exploratory LQ problems under Knightian uncertainty to quantitatively evaluate the exploration cost and its impact.
	\subsection{Non-Exploratory LQ Problems under Knightian Uncertainty}
	
	We first briefly review the infinite-horizon non-exploratory LQ control problem with discounted rewards under Knightian uncertainty. Let $B = \{ B(t),\ t \geq 0 \}$ be a one-dimensional $G$-Brownian motion defined on a sublinear expectation space. The controlled state process $\{x_t^u, t \geq 0\}$ solves the following SDE:
	\begin{equation*}
		dx_t^u = (Ax_t^u + F u_t)dt + (Cx_t^u + Du_t)dB_t, t \geq 0, x_0^u = x, 
	\end{equation*}
	where $A, F, C, D$ are given and the process $\{u_t, t \geq 0\}$ is a (non-exploratory and non-relaxed) control.
	
	The value function as defined in \eqref{eq:2}.
	
	According to \eqref{eq:10}, the corresponding HJB equation is
	\begin{align}\label{eq:22}
		\scalebox{0.9}{$\displaystyle\rho w(x)$} &\scalebox{0.9}{$\displaystyle= \frac{(2CDx\widetilde{G}[w''(x)] + Fw'(x) - Ix - Q)^2}{K - 2D^2\widetilde{G}[w''(x)]}$} \notag \\
		&\scalebox{0.9}{$\displaystyle\quad+ \left(C^2\widetilde{G}[w''(x)] - \frac{M}{2}\right)x^2 + (Aw'(x) - P)x.$} \tag{22}
	\end{align}
	
	If $K - 2D^2\widetilde{G}[w''(x)] > 0$, then the maximization on the right-hand side yields
	\[
	\boldsymbol{u}^*(x) = \frac{2CDx\widetilde{G}[w''(x)] + Fw'(x) - Ix - Q}{K - 2D^2\widetilde{G}[w''(x)]}, x \in \mathbb{R}.
	\]
	
	Through standard verification analysis, we conclude that $u$ is indeed the optimal feedback control.
	
	\subsection{Solvability Equivalence Between Non-Exploratory and Exploratory LQ Problems under Knightian Uncertainty}
	
	In this section, we prove that there exists an essential solvability equivalence between exploratory and non-exploratory LQ problems under Knightian uncertainty. Specifically, the value function and optimal control of one problem can be directly transformed into those of the other. This equivalence enables us to straightforwardly establish convergence results when the exploration weight $\lambda$ decays to zero.
	\begin{proposition}
		The following two statements (a) and (b) are equivalent.
		\begin{itemize}
			\item[(a)] The function \(V(x) = \frac{1}{2}\alpha_2 x^2 + \alpha_1 x + \alpha_0 + \frac{\lambda}{2\rho}\left(\ln\left(\frac{2\pi e \lambda}{K - \alpha_2 D^2 \overline{\sigma}^2}\right) - 1\right)\), \(x \in \mathbb{R}\), where \(\alpha_0, \alpha_1 \in \mathbb{R}\) and \(\alpha_2 < 0\), is the value function of the exploratory problem \eqref{eq:15}, and the corresponding optimal randomized feedback control is
			\begin{equation*}
				\boldsymbol{\theta}^*(u;x) = \mathcal{N} \left( u \big|\mu_3, \sigma_3^2 \right), 
			\end{equation*}
			where
			\[
			\mu_3 = \frac{(\alpha_2(F + CD\overline{\sigma}^2) - I)x + \alpha_1F - Q}{K - \alpha_2D^2\overline{\sigma}^2},\]
			\[\sigma_3^2 = \frac{\lambda}{K - \alpha_2D^2\overline{\sigma}^2}.
			\]
			
			\item[(b)] The functions \(V^{ne}(x) = \frac{1}{2}\alpha_2 x^2 + \alpha_1 x + \alpha_0\), \(x \in \mathbb{R}\), where \(\alpha_0, \alpha_1 \in \mathbb{R}\) and \(\alpha_2 < 0\), is the value function for the non-exploratory problem \eqref{eq:2}, and the corresponding optimal feedback control is
			\[
			\boldsymbol{u}^*(x) = \frac{(\alpha_2(F + CD\overline{\sigma}^2) - I)x + \alpha_1 F - Q}{K - \alpha_2 D^2 \overline{\sigma}^2}.
			\]
		\end{itemize}
	\end{proposition}
	
	\begin{proof}
		We first note that when (a) holds, the function \(v\) solves the HJB equation \eqref{eq:17} of the exploratory LQ problem. Similarly for \(w\) of the classical LQ problem when (b) holds.
		
		Next, we prove the equivalence between (a) and (b). First, a comparison between the two HJB equations \eqref{eq:17} and \eqref{eq:22} yields that if \(v\) in (a) solves the former, then \(w\) in (b) solves the latter, and vice versa.
		
		Then it remains to prove the equivalence of admissibility. Once this is established, Proposition 2 follows.
		
		The desired equivalence of the admissibility then follows from the following lemma.
		\begin{lemma}
			(i) $\liminf_{T \to \infty} e^{-\rho T} \hat{\mathbb{E}}[(X_T^*)^2] = 0$ if and only if $\liminf_{T \to \infty} e^{-\rho T} \hat{\mathbb{E}}[(x_T^*)^2] = 0$ \\
			and (ii) $\hat{\mathbb{E}}\left[\int_0^\infty e^{-\rho t} (X_t^*)^2 dt\right] < \infty$ if and only if $\hat{\mathbb{E}}\left[\int_0^\infty e^{-\rho t} (x_t^*)^2 dt\right] < \infty$.
		\end{lemma}
		The proof of Lemma 1 is provided in Appendix A. With Lemma 1 thus proved, and in conjunction with the established equivalence of the HJB equations, the proof of Proposition 2 is complete.
	\end{proof}
	It is important to note that the proof by Wang et al. \cite{wang2020reinforcement} relies on properties of linear expectation, particularly $\mathbb{E}[dX_t^2] = d(\mathbb{E}[X_t^2])$, which allows for exact solutions to the ODEs satisfied by $\mathbb{E}[X_t^2]$ and $\mathbb{E}[x_t^2]$, followed by comparison of their analytical solutions.
	
	In contrast, our proof is developed within the sublinear expectation framework, where the equality $\hat{\mathbb{E}}[dX_t^2] = d(\hat{\mathbb{E}}[X_t^2])$ does not hold, making it impossible to derive exact ODEs for $\hat{\mathbb{E}}[X_t^2]$. Consequently, we employed the $G$-Itô's formula and Gronwall's inequality to prove $\hat{\mathbb{E}}[X_T^2] \leq Y_T$ and $\hat{\mathbb{E}}[x_T^2] \leq Z_T$. Since the exponential growth rates $\alpha$ and $\tilde{\alpha}$ of these two dominating functions are identical, they guarantee the same asymptotic behavior, thereby proving the equivalence. (The specific definitions of $\alpha$, $Y_t$, $\tilde{\alpha}$ and $Z_t$ can be found in (A.1), (A.2), (A.3) and (A.4) of Appendix A.)
	
	The equivalence between statements (a) and (b) above implies that: If one problem is solvable, then the other must also be solvable; conversely, if one problem is unsolvable, the other is necessarily unsolvable.
	\subsection{Exploration Cost under Knightian Uncertainty}
	
	The exploration cost for general RL problems under Knightian uncertainty is defined as the difference between the discounted cumulative rewards under the optimal open-loop control for the non-exploratory objective \eqref{eq:2} and the exploratory objective \eqref{eq:15}, minus the entropy term. Thus, the exploration cost measures the performance loss in the non-exploratory objective due to exploration. Note that the solvability equivalence established in the previous subsection is crucial for this definition, particularly because the cost is well-defined only when both the non-exploratory and exploratory problems are solvable.
	
	Specifically, let the non-exploratory maximization problem \eqref{eq:2} have value function $V^{ne}(\cdot)$ and optimal deterministic policy $\{u_t^*, t \geq 0\}$ and let the corresponding exploratory problem \eqref{eq:15} have value function $V(\cdot)$ and optimal randomized control $\{\theta_t^*, t \geq 0\}$. Then, the exploration cost is defined as:
	\begin{equation}\label{eq:23}
		\small % 或使用 \scriptstyle 获得更小尺寸
		\begin{aligned}
			&C^{u^*, \theta^*}(x)\\
			&= \left( V^{ne}(x) - V(x) \right) \\
			&\quad - \lambda \mathcal{E} \left[ \int_0^\infty e^{-\rho t} \left( \int_U \theta_t^*(u) \ln \theta_t^*(u) du \right) dt \bigg| X_0^{\theta^*} = x \right].
		\end{aligned}
		\tag{23}
	\end{equation}
	
	Next we compute the exploration cost of the LQ problem under Knightian uncertainty. As we show, this cost is very simple: it depends only on two ``agent-specific'' parameters: the temperature parameter $\lambda$ and the discounting parameter $\rho$.
	\begin{proposition}
		Assume that statement (a) (or equivalently, (b)) of Proposition 2 holds. Then, the exploration cost of the LQ problem under Knightian uncertainty is
		\[
		\mathcal{C}^{u^*, \theta^*}(x) = \frac{\lambda}{2\rho}, \forall x \in \mathbb{R}.
		\]
	\end{proposition}
	
	\begin{proof}
		Let $\{\theta_t^*, t \geq 0\}$ be the open-loop randomized control generated by the randomized feedback control $\boldsymbol{\theta}^*$ given in statement (a) with respect to the initial state $x$, i.e.,
		\[
		\theta_t^*(u) = \mathcal{N}\left(u \big| \mu_4,\sigma_4^2 \right),
		\]
		where
		\[
		\mu_4 = \frac{(\alpha_2(F + CD\overline{\sigma}^2) - I)X_t^* + \alpha_1 F - Q}{K - \alpha_2 D^2 \overline{\sigma}^2},\]
		\[ \sigma_4^2 = \frac{\lambda}{K - \alpha_2 D^2 \overline{\sigma}^2}.
		\]
		Here, $\{X_t^*, t \geq 0\}$ is the state process starting from state $x$ under the application of $\boldsymbol{\theta}^*$ for the exploratory problem. It then straightforwardly follows that
		\[
		\begin{aligned}
			\int_{\mathbb{R}} \theta_t^*(u) \ln \theta_t^*(u) du 
			&= -\frac{1}{2}\left( \ln(2\pi\sigma_4^2) + 1 \right).
		\end{aligned}
		\]
		Substituting $\sigma_4^2$ into the above equation, we get
		\[
		\int_{\mathbb{R}} \theta_t^*(u) \ln \theta_t^*(u) du = -\frac{1}{2} \ln \left( \frac{2\pi e \lambda}{K - \alpha_2 D^2 \overline{\sigma}^2} \right).
		\]
		
		The desired result now follows directly from the general definition in \eqref{eq:23} and the expressions in Proposition 2 (a) for $V(\cdot)$ and (b) for $V^{ne}(\cdot)$. 
		Thus, the proof is complete.
	\end{proof}
	
	The exploration cost under Knightian uncertainty obtained in this paper seems counterintuitive, yet it is consistent with the result derived by Wang et al. \cite{wang2020reinforcement} under the classical linear expectation framework. The reason is that, when calculating the cost, the expectation of the entropy term and the discrepancy of the value function exactly cancel out all system-dependent components, leaving only constants related to $\lambda$ and $\rho$. Therefore, even though the uncertainty faced by the agent is elevated from ordinary risk to the more fundamental level of Knightian uncertainty, the marginal cost of exploration is still determined solely by the agent's own preference parameters—the exploration weight $\lambda$ and the discount rate $\rho$—and is independent of the degree of external environmental uncertainty.
	The finding provides stronger theoretical justification for employing randomized policies in complex, uncertain environments.
	
	\subsection{Vanishing Exploration}
	The exploration weight $\lambda$ is considered as an exogenous parameter reflecting the level of exploration expected by the RL agent. The smaller this parameter is, the more emphasis is placed on exploitation. When this parameter is close enough to zero, the exploratory formulation is close enough to the problem without exploration. Naturally, the ideal result is that if the exploration weight $\lambda$ decays zero, then the exploratory LQ problem under Knightian uncertainty will converge to its non-exploratory counterpart. The following result makes this precise.
	\begin{proposition}
		Assumes that statement (a) (or equivalently, (b)) of Proposition 2 holds. For each $x \in \mathbb{R}$,
		\[
		\lim_{\lambda \to 0} \boldsymbol{\theta}^*(\cdot; x) = \delta_{\boldsymbol{u}^*(x)}(\cdot) \quad \text{weakly}.
		\]
		In addition, for each $x \in \mathbb{R}$,
		\[
		\lim_{\lambda \to 0} |V(x) - V^{ne}(x)| = 0.
		\]
	\end{proposition}
	\begin{proof}
		According to Proposition 2, we have $\mu_3 = \boldsymbol{u}^*(x)$, which means the mean of the randomized strategy exactly equals the deterministic optimal control. And the variance of the randomized strategy is
		\[\sigma_3^2 = \frac{\lambda}{K - \alpha_2D^2\overline{\sigma}^2}.
		\]
		So we get
		\begin{equation*}
			\lim_{\lambda \to 0} \sigma_3^2 = 0.
		\end{equation*}
		For a Gaussian distribution $\mathcal{N}(\mu_3, \sigma_3^2)$, when the variance $\sigma_3^2 \to 0$, it converges weakly to the Dirac delta function $\delta_u$. Therefore, 
		\[
		\lim_{\lambda \to 0} \boldsymbol{\theta}^*(\cdot; x) = \delta_{\boldsymbol{u}^*(x)}(\cdot) \quad \text{weakly}.
		\]
		According to Proposition 2, we have
		\begin{align*}
			V(x) - V^{\text{ne}}(x) = \frac{\lambda}{2\rho} \left( \ln \left( \frac{2\pi e \lambda}{K - \alpha_2 D^2 \overline{\sigma}^2} \right) - 1 \right).
		\end{align*}
		Since
		\[
		\lim_{\lambda \to 0} \frac{\lambda}{2\rho} \left( \ln \left( \frac{2\pi e}{K - \alpha_2 D^2 \overline{\sigma}^2} \right) - 1 \right) = 0,
		\]
		we get
		\[
		\lim_{\lambda \to 0} \left| V(x) - V^{\text{ne}}(x) \right| = 0.
		\]
		Thus, Proposition 4 is proved.
	\end{proof}
	We wish to emphasize that the proof of this limit result is not a straightforward application of existing standard tools. Its validity critically depends on the non-trivial results established within our sublinear expectation framework, specifically Proposition 2. The key inputs---the equality $\mu_3 = \boldsymbol{u^*(x)}$ and the form of $\sigma_3^2$---both are results derived within the sublinear expectation framework in this paper. The convergence of the value function also relies on the explicit expressions of the exploratory value function and non-exploratory value function provided in Proposition 2. These expressions themselves are main results of our work, obtained by using the new tools of sublinear expectation, and are not attainable via classical linear expectation tools (e.g., the methods in \cite{wang2020reinforcement}).
	
In Appendix C, we provide a practical LQ example to verify the theoretical results of this paper via numerical simulations, with the verification framework for the simulations detailed in Appendix B.
	
	\section{Conclusion}
	This paper establishes a unified theoretical framework for studying optimal control via RL under Knightian uncertainty by integrating sublinear expectation theory with entropy-regularized relaxed stochastic control. Our work fundamentally extends the research paradigm based on linear expectation established by Wang et al. \cite{wang2020reinforcement}, providing a more solid theoretical foundation for decision-making in environments with Knightian uncertainty.
	
	Although this paper is primarily theoretical and does not propose a new learning algorithm, our findings provide a foundational roadmap for developing practical RL algorithms under Knightian uncertainty. The theoretical results derived here can directly inform algorithm design in the following ways:
	
	\begin{itemize}
		\item System Identification: The volatility bounds $[\underline{\sigma}^2, \overline{\sigma}^2]$ are key parameters for quantifying Knightian uncertainty. The algorithm needs to incorporate a parallel online estimation module that employs the $\varphi$-max-mean algorithm proposed by Peng \cite{peng2017} to estimate these bounds in real-time, and use them as regulatory signals input to the networks.
		
		\item Policy Evaluation: To meet the decision-making requirements of an ambiguity-averse agent, the paradigm of policy evaluation must shift from traditional ``performance evaluation'' to ``robustness evaluation.'' Specifically, the objective is no longer to assess the expected value of a policy in the average environment, but rather its lower expectation value under the worst-case scenario. This necessitates a reformulation of the update rule for the Critic. For instance, in temporal difference learning, the target value should be altered from $\mathbb{E}[r + \gamma V(s')]$ to $\mathcal{E}[r + \gamma V(s')]$, thereby providing a robust value baseline for policy improvement.
		
		\item Policy Improvement: The actor network must output a complete Gaussian policy, where both the mean and variance are dynamically modulated according to environmental uncertainty. Specifically, the mean calculation must incorporate uncertainty parameters to output the most robust action, while the variance serves as an adaptive exploration parameter that is negatively correlated with the level of uncertainty. This achieves an endogenous risk-aware exploration mechanism that surpasses traditional fixed entropy bonus-based exploration methods.

	\end{itemize}
	
	These specific design principles will form the basis for developing practical algorithms in our immediate subsequent work.
	
	Our future research will pursue two directions. First, we will apply the exploratory state dynamics under Knightian uncertainty developed in this paper to continuous-time mean-variance portfolio selection. By casting this problem as a trade-off between exploration and exploitation, we aim to derive optimal investment strategies for Knightian uncertain environments. Second, we will leverage real financial market data to online estimate the bounds of volatility and develop novel reinforcement learning algorithms to quantitatively analyze the specific impact mechanisms of Knightian uncertainty on investment decisions.
	
	\section*{APPENDIX}

	% 附录A：引理1证明
	\section*{Appendix A: Proof of Lemma 1}
	\begin{proof}
		Throughout this proof, we let \(x\) be fixed, being the initial state of both the exploratory problem in statement (a) and the classical problem in statement (b). Let \(\theta^{*}=\{\theta^{*}_{t},t\geq 0\}\) and \(u^{*}=\{u^{*}_{t},t\geq 0\}\) be respectively the open-loop controls generated by the feedback controls \(\boldsymbol{\theta}^{*}\) and \(\boldsymbol{u}^{*}\) of the two problems, and \(X^{*}=\{X^{*}_{t},t\geq 0\}\) and \(x^{*}=\{x^{*}_{t},t\geq 0\}\) be respectively the corresponding state processes, both starting from \(x\).
		
		To ease the presentation, we rewrite the exploratory dynamics of \(X^{*}\) under \(\theta^{*}\) as:
		\begin{equation*}
			\begin{aligned}
				dX_t^* &= \bigg[ \Big( A + \tfrac{F(\alpha_2(F + CD\overline{\sigma}^2) - I)}{K - \alpha_2D^2\overline{\sigma}^2} \Big) X_t^* + \tfrac{F(\alpha_1F - Q)}{K - \alpha_2D^2\overline{\sigma}^2} \bigg] dt \\
				&\quad+ \bigg[ \Big( \big( C + \tfrac{D(\alpha_2(F + CD\overline{\sigma}^2) - I)}{K - \alpha_2D^2\overline{\sigma}^2} \big) X_t^* + \tfrac{D(\alpha_1F - Q)}{K - \alpha_2D^2\overline{\sigma}^2} \Big)^2 \\
				&\quad+ \tfrac{\lambda D^2}{K - \alpha_2D^2\overline{\sigma}^2} \bigg]^{1/2} dB_t\\
				&=(A_{1}X_{t}^{*}+A_{2})dt+\sqrt{(\mathcal{B}_{1}X_{t}^{*}+\mathcal{B}_{2})^{2}+C_{1}}dB_{t},
			\end{aligned}
		\end{equation*}
		where
		\[
		\begin{aligned}
			A_{1} &= A+\frac{F(\alpha_{2}(F+CD\overline{\sigma}^2)-I)}{K-\alpha_{2}D^{2}\overline{\sigma}^2}, \\
			A_{2} &:= \frac{F(\alpha_{1}F-Q)}{K-\alpha_{2}D^{2}\overline{\sigma}^2}, \\
			\mathcal{B}_{1} &:= C+\frac{D(\alpha_{2}(F+CD\overline{\sigma}^2)-I)}{K-\alpha_{2}D^{2}\overline{\sigma}^2}, \\
			\mathcal{B}_{2} &:= \frac{D(\alpha_{1}F-Q)}{K-\alpha_{2}D^{2}\overline{\sigma}^2}, \\
			C_{1} &:= \frac{\lambda D^{2}}{K-\alpha_{2}D^{2}\overline{\sigma}^2}.
		\end{aligned}
		\]
		Similarly, the classical dynamics of \(x^{*}\) under \(u^{*}\) solves
		\[
		dx_{t}^{*}=(A_{1}x_{t}^{*}+A_{2})dt+(\mathcal{B}_{1}x_{t}^{*}+\mathcal{B}_{2})dB_{t}.
		\]
		
		Let $f(x) = x^2$. Applying the \textit{G}-Itô's formula to the exploratory process $X_t$, we get
		\[
		d(X_t^2) = 2X_t dX_t + d\langle X \rangle_t.
		\]
		Recall the dynamics
		\[
		\begin{aligned}
			dX_t &= (A_1 X_t + A_2)dt + \sigma_t dB_t, \\
			\sigma_t &= \sqrt{(\mathcal{B}_1 X_t + \mathcal{B}_2)^2 + C_1}.
		\end{aligned}
		\]
		Then, we have
		\[
		\begin{aligned}
			d(X_t^2) &= \left[ 2X_t(A_1 X_t + A_2) + \sigma_t^2 \right] dt + 2X_t \sigma_t \, dB_t \\
			&\quad + \sigma_t^2 \, d\langle B \rangle_t.
		\end{aligned}
		\]
		Since $d\langle B \rangle_t \leq \overline{\sigma}^2 dt$ holds quasi-surely (q.s.), we obtain
		\[
		d(X_t^2) \leq \left[2A_1 X_t^2 + 2A_2 X_t + \sigma_t^2  \overline{\sigma}^2\right] dt + 2X_t \sigma_t\,dB_t.
		\]
		So
		\begin{align*}
			d(X_t^2) &\leq \left[ (2A_1 + \overline{\sigma}^2 \mathcal{B}_1^2) X_t^2 + (2A_2 + 2\overline{\sigma}^2 \mathcal{B}_1 \mathcal{B}_2) X_t \right. \\
			&\quad \left. + \overline{\sigma}^2(\mathcal{B}_2^2 + C_1) \right] dt + 2X_t \sigma_t dB_t.
		\end{align*}
		Let
		\begin{align}
			\alpha &= 2A_1 + \overline{\sigma}^2 \mathcal{B}_1^2 + \frac{|2A_2 + 2\overline{\sigma}^2 \mathcal{B}_1 \mathcal{B}_2|}{2}, \label{A.1}\tag{A.1}\\
			\beta  &= \frac{|2A_2 + 2\overline{\sigma}^2 \mathcal{B}_1 \mathcal{B}_2|}{2} +\overline{\sigma}^2 (\mathcal{B}_2^2 + C_1). \nonumber
		\end{align}
		Then, we have
		\[
		d(X_t^2) \leq [\alpha X_t^2 + \beta]\, dt + 2X_t \sigma_t\, dB_t.
		\]
		Integrating both sides from 0 to $t$, we have
		\[
		\int_0^t d(X_s^2) \leq \int_0^t [\alpha X_s^2 + \beta]\, ds + \int_0^t 2X_s \sigma_s\, dB_s.
		\]
		Then, we get
		\[
		X_t^2 \leq x^2 + \int_0^t (\alpha X_s^2 + \beta)  ds + \int_0^t 2X_s\sigma_s  dB_s.
		\]
		Taking expectation $\mathbb{\hat{E}}$ on both sides, we deduce
		\[
		\hat{\mathbb{E}}[X_t^2] \leq x^2 + \hat{\mathbb{E}}\left[ \int_0^t (\alpha X_s^2 + \beta) ds \right].
		\]
		Denoting $m(t) = \hat{\mathbb{E}}[X_t^2]$, we have
		\begin{align*}
			m(t) &\leq x^2 + \hat{\mathbb{E}} \left[ \int_0^t (\alpha X_s^2 + \beta)  ds \right], \\
			m(t) &\leq x^2 + \int_0^t (\alpha m(s) + \beta)  ds.
		\end{align*}
		Thus, we get
		\[
		m(t) \leq x^2 + \beta t + \alpha \int_0^t m(s)  ds.
		\]
		Applying Gronwall's inequality to the above integral inequality, we can obtain
		\[
		m(t) \leq e^{\alpha t} x^{2} + \frac{\beta}{\alpha} (e^{\alpha t} - 1).
		\]
		Denoting
		\begin{equation*}
			Y_t = e^{\alpha t} x^2 + \frac{\beta}{\alpha}(e^{\alpha t} - 1)\label{A.2}\tag{A.2},
		\end{equation*}
		we have
		\[
		\hat{\mathbb{E}}[X_T^2] \leq Y_T.
		\]
		Therefore, we have
		\[
		\begin{aligned}
			e^{-\rho T} \hat{\mathbb{E}}[X_T^2] 
			& \leq e^{-\rho T} Y_T \\
			& = e^{-(\rho - \alpha)T} \left(x^2 + \frac{\beta}{\alpha}\right) - \frac{\beta}{\alpha} e^{-\rho T}.
		\end{aligned}
		\]
		If $\rho > \alpha$, then the right-hand side tends to 0 as $T \to \infty$, and then we have
		\[
		\liminf_{T \to \infty} e^{-\rho T} \hat{\mathbb{E}}[X_T^2] = 0.
		\]
		
		For $x_t$, we perform the same analysis. The only difference is that $C_1 = 0$. Thus, we get
		\[
		d(x_t^2) \leq \left[ \tilde{\alpha} x_t^2 + \tilde{\beta} \right] dt + 2 x_t ( \mathcal{B}_1 x_t + \mathcal{B}_2 ) d B_t,
		\]
		where
		\begin{align}
			\tilde{\alpha} &= 2A_1 + \overline{\sigma}^2 \mathcal{B}_1^2 + \frac{|2A_2 + 2\overline{\sigma}^2 \mathcal{B}_1 \mathcal{B}_2|}{2}\label{A.3}, \tag{A.3} \\
			\tilde{\beta}  &= \frac{|2A_2 + 2\overline{\sigma}^2 \mathcal{B}_1 \mathcal{B}_2|}{2} + \overline{\sigma}^2 \mathcal{B}_2^2. \nonumber
		\end{align}
		Denoting $\tilde{m}(t) = \hat{\mathbb{E}}[x_t^2]$, we have
		\[
		\tilde{m}(t) \leq x^2 + \tilde{\beta} t + \tilde{\alpha} \int_0^t \tilde{m}(s)  ds.
		\]
		Applying Gronwall's inequality to the above integral inequality, we can obtain
		\[
		\tilde{m}(t) \leq e^{\tilde{\alpha} t} x^{2} + \frac{\tilde{\beta}}{\tilde{\alpha}} (e^{\tilde{\alpha} t} - 1).
		\]
		Denoting
		\begin{equation*}
			Z_t = e^{\tilde{\alpha}t} x^2 + \frac{\tilde{\beta}}{\tilde{\alpha}} (e^{\tilde{\alpha}t} - 1) \label{A.4}\tag{A.4},
		\end{equation*}
		we have
		\[
		\hat{\mathbb{E}}[x_T^2] \leq Z_T.
		\]
		Therefore, we get
		\[
		\begin{aligned}
			e^{-\rho T} \hat{\mathbb{E}}[x_T^2] 
			& \leq e^{-\rho T} Z_T \\
			& = e^{-(\rho - \tilde{\alpha})T} \left(x^2 + \frac{\tilde{\beta}}{\tilde{\alpha}}\right) - \frac{\tilde{\beta}}{\tilde{\alpha}} e^{-\rho T}.
		\end{aligned}
		\]
		If \( \rho > \tilde{\alpha} \), then the right-hand side tends to 0 as $T \to \infty$, and then we get
		\[
		\liminf_{T \to \infty} e^{-\rho T} \hat{\mathbb{E}}[x_T^2] = 0.
		\]
		Based on the above discussion, we have
		\[
		\alpha = \tilde{\alpha}.
		\]
		Hence, we get
		\[
		\rho > \alpha \quad \text{if and only if} \quad \rho > \tilde{\alpha}.
		\]
		This implies $\liminf_{T \to \infty} e^{-\rho T} \hat{\mathbb{E}}[X_T^2] = 0 \text{ if and only if }$ \\
		$\liminf_{T \to \infty} e^{-\rho T} \hat{\mathbb{E}}[x_T^2] = 0.$
		
		Moreover, we have
		\begin{align*}
			&\hat{\mathbb{E}}\left[\int_{0}^{\infty}e^{-\rho t}X_{t}^{2}dt\right] \\
			&\leq \int_{0}^{\infty}e^{-\rho t}\hat{\mathbb{E}}[X_{t}^{2}]dt \\
			&\leq \int_{0}^{\infty}e^{-\rho t}Y_{t}dt \\
			&= \int_{0}^{\infty}e^{-\rho t}\left(e^{\alpha t}x^{2}+\frac{\beta}{\alpha}(e^{\alpha t}-1)\right)dt.
		\end{align*}
		This integral converges if and only if $\rho > \alpha$.
		
		For $x_t$ the analysis is completely analogous: the integral converges if and only if $\rho > \tilde{\alpha}$.
		
		Since $\rho > \alpha \iff \rho > \tilde{\alpha}$, this proves the equivalence of the integrability conditions. Thus, $\hat{\mathbb{E}}\left[\int_0^\infty e^{-\rho t} (X_t^*)^2 dt\right] < \infty$ if and only if $\hat{\mathbb{E}}\left[\int_0^\infty e^{-\rho t} (x_t^*)^2 dt\right] < \infty$. 
		
		Therefore, Lemma 1 is proved.
	\end{proof}

	\newpage
	\section*{Appendix B: Numerical Verification Framework}
	\begin{table}[h]
		\centering
		\caption{Systematic Verification Framework}
		\label{tab:verification-framework}
		\footnotesize
		\begin{tabular}{@{}p{0.95\linewidth}@{}}
			\toprule
			\textbf{Model Parameter Definition}\\
			System parameters: $A, F, C, D, M, I, K, P, Q$\\
			Test parameters: $x_{\mathrm{test}}, \overline{\sigma}_{\mathrm{range}}, \lambda_{\mathrm{range}}, \rho_{\mathrm{range}}$\\
			Numerical settings: $\epsilon, N$\\
			\midrule
			\textbf{For each analysis mode}:\\
			\quad \textbf{Mode A: Gaussian Policy Verification and Normality Test}\\
			\quad \quad \textbf{for each} $\sigma$ in $\overline{\sigma}_{\mathrm{range}}$:\\
			\quad \quad \quad Compute $\mu, \sigma_{\mathrm{pol}} \leftarrow \mathrm{optimal\_policy}(x_{\mathrm{test}}, \overline{\sigma})$\\
			\quad \quad \quad Compute $\alpha \leftarrow \mathrm{compute\_stability\_coefficient}(\overline{\sigma})$\\
			\quad \quad \quad Verify $\rho > \alpha$\\
			\quad \quad \quad \textbf{Normality Test}:\\
			\quad \quad \quad \quad Generate $N$ samples $u_i \sim \mathcal{N}(\mu, \sigma_{\mathrm{pol}}^2)$\\
			\quad \quad \quad \quad Compute K-S and A-D tests\\
			\quad \quad \textbf{end for}\\
			\quad \quad Output: Probability distributions + Normality tests + Stability report\\[2pt]
			
			\quad \textbf{Mode B: Discount Rate Impact Analysis}\\
			\quad \quad \textbf{for each} $\rho$ in $\rho_{\mathrm{range}}$:\\
			\quad \quad \quad Compute policy parameters\\
			\quad \quad \quad Verify stability conditions\\
			\quad \quad \textbf{end for}\\
			\quad \quad Analyze: $\partial \mathrm{Var}/\partial \rho < 0$\\
			\quad \quad Output: Variance sensitivity plots + Analysis report\\[2pt]
			
			\quad \textbf{Mode C: Distribution Convergence ($\lambda \to 0$)}\\
			\quad \quad \textbf{for each} $\sigma$ in $\overline{\sigma}_{\mathrm{range}}$:\\
			\quad \quad \quad \textbf{for each} $\lambda$ in $\lambda_{\mathrm{range}}$:\\
			\quad \quad \quad \quad Compute policy parameters\\
			\quad \quad \quad \quad Track $\sigma_{\mathrm{pol}}(\lambda)$\\
			\quad \quad \quad \textbf{end for}\\
			\quad \quad \quad Verify $\sigma_{\mathrm{pol}}(\lambda)$ decreases monotonically as $\lambda \to 0$\\
			\quad \quad \textbf{end for}\\
			\quad \quad Output: Convergence plots + Distribution evolution\\[2pt]
			
			\quad \textbf{Mode D: Value Function Convergence ($\lambda \to 0$)}\\
			\quad \quad \textbf{for each} $\sigma$ in $\overline{\sigma}_{\mathrm{range}}$:\\
			\quad \quad \quad \textbf{for each} $\lambda$ in $\lambda_{\mathrm{range}}$:\\
			\quad \quad \quad \quad Compute $V(x)$ and $V^{\mathrm{ne}}(x)$\\
			\quad \quad \textbf{end for}\\
			\quad \quad \quad Verify $V(x)$ approaches $V^{\mathrm{ne}}(x)$ as $\lambda \to 0$\\
			\quad \quad \textbf{end for}\\
			\quad \quad Output: Value function convergence plots\\
			\bottomrule
		\end{tabular}
	\end{table}
	
	% 第二部分：核心函数
	\begin{table}[H]
		\centering
		\caption{Computational Components}
		\label{tab:computational-components}
		\footnotesize
		\begin{tabular}{@{}p{0.95\linewidth}@{}}
			\toprule
			\textbf{Computational Components}\\
			\midrule
			1. $\mathrm{solve\_hjb}(\overline{\sigma}, \rho)$:\\
			\quad Solve $f(k_2) = \rho k_2 - \dfrac{2(CD \cdot k_2 - I)^2}{\rho(K - k_2 D^2 \overline\sigma^2)} - \rho[(C^2\overline\sigma^2 + 2A)k_2 - M]$\\
			\quad $k_1 \leftarrow \dfrac{-2(CD k_2 - I)Q - \rho P}{\rho(K - k_2 D^2 \overline\sigma^2) - 2F(CD k_2 - I) - \rho A}$\\
			\quad $k_0 \leftarrow -\dfrac{Q^2 + \lambda \ln(2\pi e\lambda/(K - k_2 D^2 \overline\sigma^2))}{2\rho(K - k_2 D^2 \overline\sigma^2)}$\\[4pt]
			
			2. $\mathrm{compute\_stability\_coefficient}(\overline{\sigma})$:\\
			\quad Compute $\alpha$ as defined in (B.1) of Appendix B\\
			\quad $\alpha = 2A_1 + \overline{\sigma}^2 \mathcal{B}_1^2 + \frac{|2A_2 + 2\overline{\sigma}^2 \mathcal{B}_1 \mathcal{B}_2|}{2}$\\[4pt]
			
			3. $\mathrm{optimal\_policy}(x_{\mathrm{test}}, \overline{\sigma})$:\\
			\quad $\mu \leftarrow \dfrac{k_2(F + CD\overline\sigma^2) - I}{K - k_2 D^2 \overline\sigma^2}x + \dfrac{k_1 F - Q}{K - k_2 D^2 \overline\sigma^2}$\\
			\quad $\sigma_{\mathrm{pol}} \leftarrow \sqrt{\dfrac{\lambda}{K - k_2 D^2 \overline\sigma^2}}$\\[4pt]
			
			4. $\mathrm{exploratory\_value\_function}(x_{\mathrm{test}}, \overline{\sigma})$:\\
			\quad $V \leftarrow \dfrac{1}{2}k_2 x^2 + k_1 x + k_0 + \dfrac{\lambda}{2\rho}\left[\ln\dfrac{2\pi e\lambda}{K - k_2 D^2 \overline\sigma^2} - 1\right]$\\[4pt]
			
			5. $\mathrm{non-exploratory\_value\_function}(x_{\mathrm{test}})$:\\
			\quad $V \leftarrow \dfrac{1}{2}k_2 x^2 + k_1 x + k_0$ (with $\lambda \to 0$)\\
			\bottomrule
		\end{tabular}
	\end{table}
	
	\newpage
	% 第三部分：验证指标和输出
	\begin{table}[h!]
		\centering
		\caption{Verification Criteria and Outputs}
		\label{tab:verification-criteria}
		\footnotesize
		\begin{tabular}{@{}p{0.95\linewidth}@{}}
			\toprule
			\textbf{Verification Criteria}\\
			\midrule
			\textbf{Normality}: K-S test $p$-value $> 0.05$, A-D statistic $<$ critical value\\
			\textbf{Stability}: $\rho > \alpha$, $k_2 < 0$, $K - k_2 D^2 \overline\sigma^2 > 0$\\
			\textbf{Convergence}: $\lim_{\lambda \to 0} \sigma_{\text{pol}} = 0$, $\lim_{\lambda \to 0} V(x) = V^{ne}(x)$\\
			\textbf{Sensitivity}: $\dfrac{\partial \mathrm{Var}}{\partial \rho} < 0$, $\dfrac{\partial \mathrm{Var}}{\partial \overline{\sigma}} < 0$\\
			\midrule
			\textbf{Output}\\
			\midrule
			Policy distribution visualizations, parameter sensitivity analysis plots, \\
			convergence verification plots, statistical test results, and stability assessment reports\\\\
			\bottomrule
		\end{tabular}
	\end{table}
	
	\section*{Appendix C: Numerical Simulation}
	To verify the preceding theoretical results, this Appendix presents a simple LQ example where the reward depends on both the state and the action, with numerical simulations implemented using Jupyter Notebook. Please refer to Appendix B for the complete numerical verification framework.
	
	\subsection{Experimental Setup and Parameter Configuration}
	Consider an indoor temperature control system where the state variable $x$ represents the deviation of the indoor temperature from the setpoint (unit: $^\circ$C). The control variable $u$ denotes the heater's power (unit: kW), and the system is subject to random disturbances from ambient temperature fluctuations and device noise. The system state dynamic is simulated according to (1) and (12), while the reward function is simulated based on (13). The control objective is to achieve precise regulation of indoor temperature while balancing energy efficiency, stability, comfort, and cost-effectiveness, thereby providing occupants with a comfortable, energy-efficient, and stable indoor environment.
	
	The symbols and definitions of the specific parameters are detailed in Table \ref{tab:settings}.
	\begin{table}[!h] 
		\caption{Basic Parameter Settings}
		\label{tab:settings}
		\centering
		\renewcommand{\arraystretch}{1.1} 
		\setlength{\tabcolsep}{4pt} 
		\begin{tabularx}{\columnwidth}{@{}>{\raggedright\arraybackslash}X 
				>{\centering\arraybackslash}m{1.8cm}
				>{\centering\arraybackslash}m{1.5cm}@{}}
			\toprule
			\textbf{Parameter Description} & \textbf{Symbol} & \textbf{Value} \\ 
			\midrule
			Natural decay rate of indoor temperature & $A$ & -0.2 \\ 
			Heater efficiency & $F$ & 0.8 \\
			State-dependent noise & $C$ & 0.5 \\
			Control-dependent noise & $D$ & 1.2 \\
			\addlinespace[0.1cm]
			
			State quadratic penalty & $M$ & 10 \\
			Control quadratic penalty & $K$ & 2 \\
			State linear penalty & $P$ & 0.5 \\
			Control linear penalty & $Q$ & 0.2 \\
			Cross-term penalty & $I$ & 0.3 \\
			\addlinespace[0.1cm]
			
			Number of samples & $N$ & 10,000 \\ 
			\bottomrule
		\end{tabularx}
	\end{table}
	
	Substituting the above parameters into (1) and (12), the system dynamics is governed by
	\[
	dx_t^u = (-0.2x_t^u + 0.8u_t)dt + (0.5x_t^u + 1.2u_t)dB_t, t \geq 0, 
	\]
	where $x_0^u = x$.
	
	Substituting the above parameters into (13), the reward function as
	\[
	r(x, u) = -\left(5x^2 + 0.3xu + u^2 + 0.5x + 0.2u\right),
	\]
	where $x, u \in \mathbb{R}$.
	
\begin{figure}[!h]
		\centering
		\includegraphics[width=0.5\linewidth]{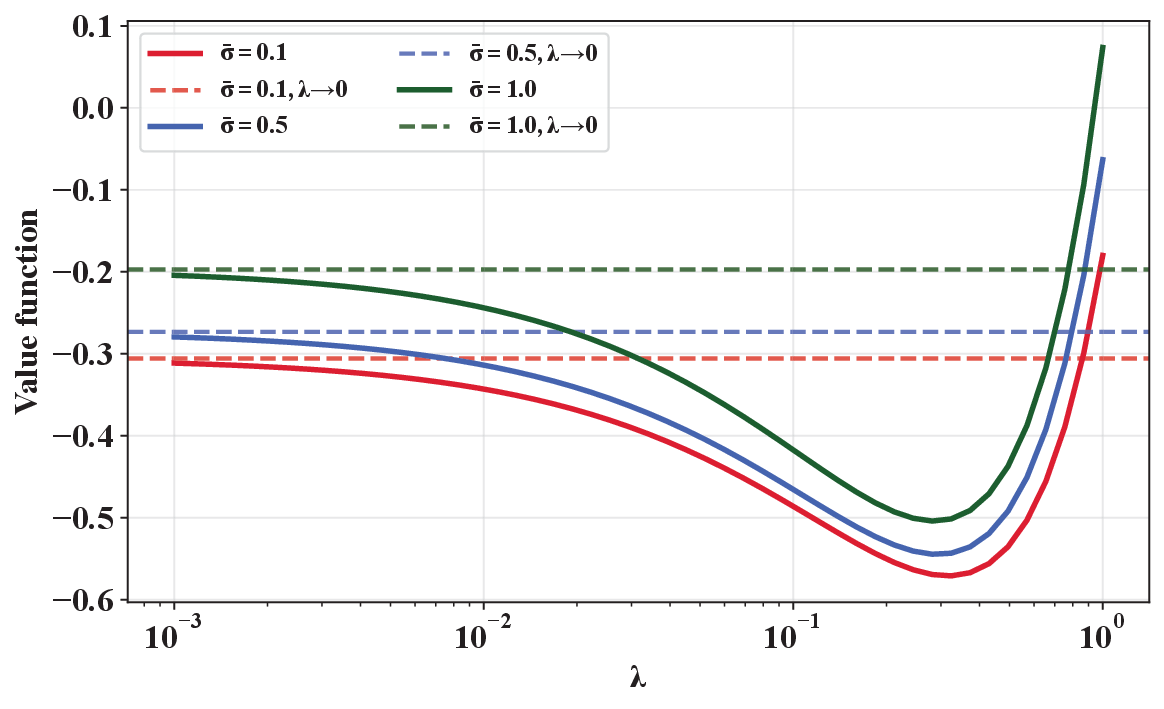} % 单栏适配：宽度改为0.8\linewidth
		\caption{Convergence of the exploratory value function $v(x)$}
		\label{fig:5}
	\end{figure}

	As established in Proposition 2, the solvability equivalence between exploratory and non-exploratory problems necessitates the discount rate $\rho$ to satisfy the conditions $\rho > \alpha$ (The specific definition of $\alpha$ can be found in (A.1) of Appendix A). Moreover, since $\rho$ indirectly affects $k_2$, we fix $\rho = 0.3$ in Figs. \ref{fig:1}, \ref{fig:2}, \ref{fig:4}, and \ref{fig:5} to simultaneously satisfy the solvability condition $\rho > \alpha$ (where $\alpha < 0$ under our parameter configuration) and the concavity requirement of the value function ($k_2 < 0$). Thus, all subsequent simulations are  conducted under the validated condition $k_2 < 0$, with a fixed temperature deviation of $x = 1^\circ$C to examine policy behavior around the setpoint.
	
	\subsection{Validation of Optimal Policy Distribution Characteristics}
\begin{figure}[H]
		\centering
		\includegraphics[width=0.5\linewidth]{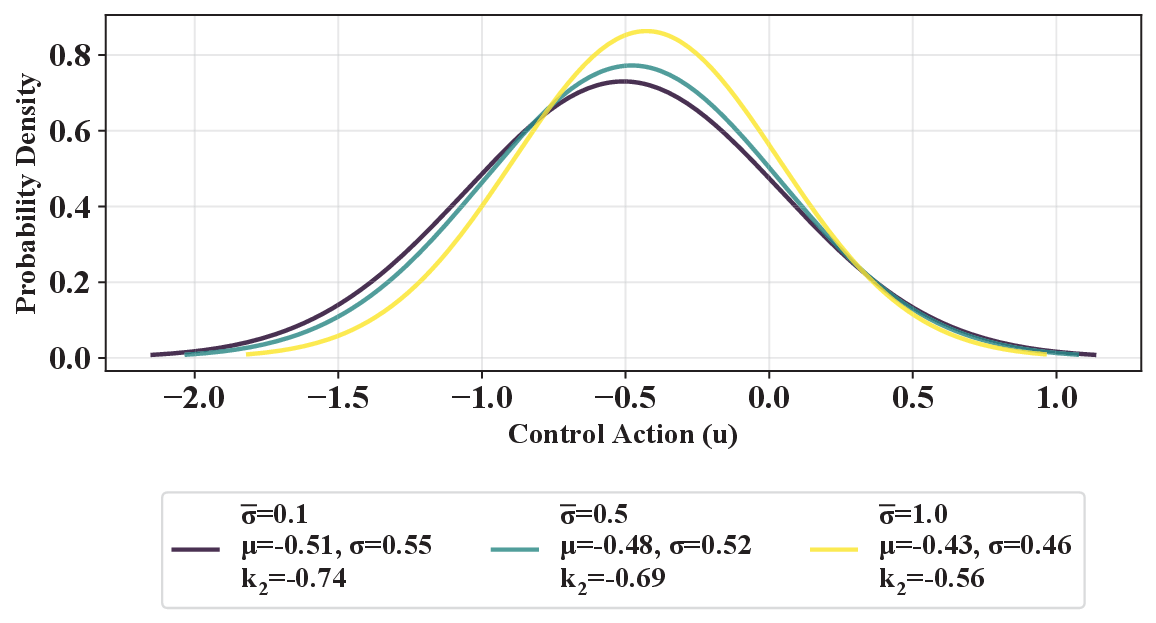} % 单栏适配：宽度改为0.8\linewidth
		\caption{Optimal policy distribution}
		\label{fig:1}
	\end{figure}

\begin{figure}[!h]
		\centering
		\includegraphics[width=0.5\linewidth]{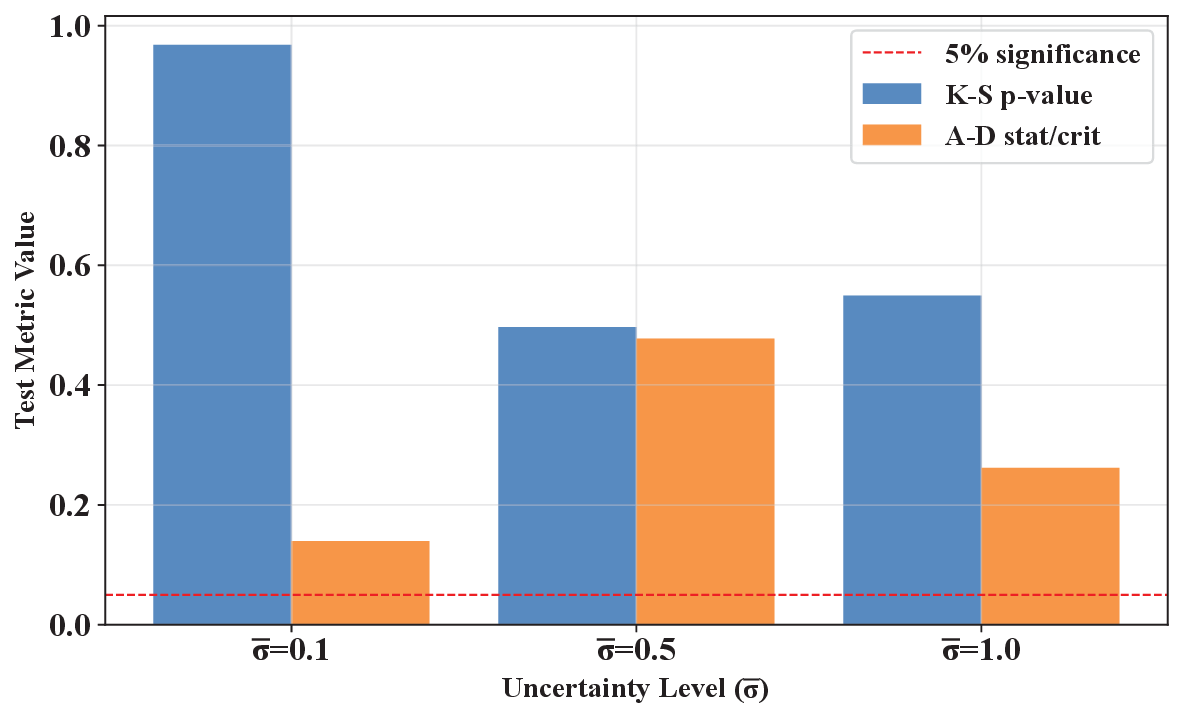} % 单栏适配：宽度改为0.8\linewidth
		\caption{Normality test results}
		\label{fig:2}
	\end{figure}

	In Figs. \ref{fig:1} and~\ref{fig:2}, with $\lambda$ fixed at 0.6, we analyze the distribution characteristics of the optimal control policy for this indoor temperature regulation system under different $\overline{\sigma}$ values (0.1, 0.5, 1.0).
	
	As shown in Fig.~\ref{fig:1}, the three colored curves represent the optimal policy distributions under different $\overline{\sigma}$ values. Each curve exhibits a bell shape, indicating that the optimal policies under varying Knightian uncertainty levels appear to follow Gaussian distributions.
	
	Fig.~\ref{fig:2} presents the normality test results: 
	\begin{itemize}
		\item The blue bars indicate the p-values from the Kolmogorov-Smirnov (K-S) tests for different $\overline{\sigma}$ values. All p-values are significantly above the 0.05 threshold (red line), failing to reject the null hypothesis that the optimal policy follow Gaussian distribution.
		
		\item The orange bars represent the Anderson-Darling (A-D) test (statistic-to-critical-value ratios). All ratios are below 1, showing that the A-D statistics are smaller than their respective critical values, further supporting the normality assumption.
	\end{itemize}
	
	Integrating the results from Figs.~\ref{fig:1}--\ref{fig:2}, we conclude that within the LQ framework, the numerically simulated optimal policy exhibit perfect agreement with theoretical prediction. This validates both the formulations in~(18)--(21) and Proposition 1, demonstrating that the optimal randomized feedback control policy follows a Gaussian distribution across all tested Knightian uncertainty levels under state-dependent reward structures.
	
	Furthermore, we observe that in Fig.~\ref{fig:1}, the purple curve exhibits the widest distribution, indicating the greatest data dispersion (largest variance), while the yellow curve appears the most concentrated, reflecting the most compact data distribution (smallest variance). This phenomenon demonstrates that as $\overline{\sigma}$ increases (corresponding to higher Knightian uncertainty), the variance decreases and the exploration range narrows. Consequently, the agent of the temperature control system adopts a more conservative strategy by prioritizing current optimal control actions over extensive exploration, thereby achieving a balance between cost-efficiency and stability. These findings are fully consistent with our theoretical analysis based on the LQ framework with state-dependent rewards, confirming the modulating effect of Knightian uncertainty on the variance of optimal feedback policies.
	
	\subsection{The Impact of $\rho$ on the Optimal Policy}
\begin{figure}[H]
		\centering
		\includegraphics[width=0.5\linewidth]{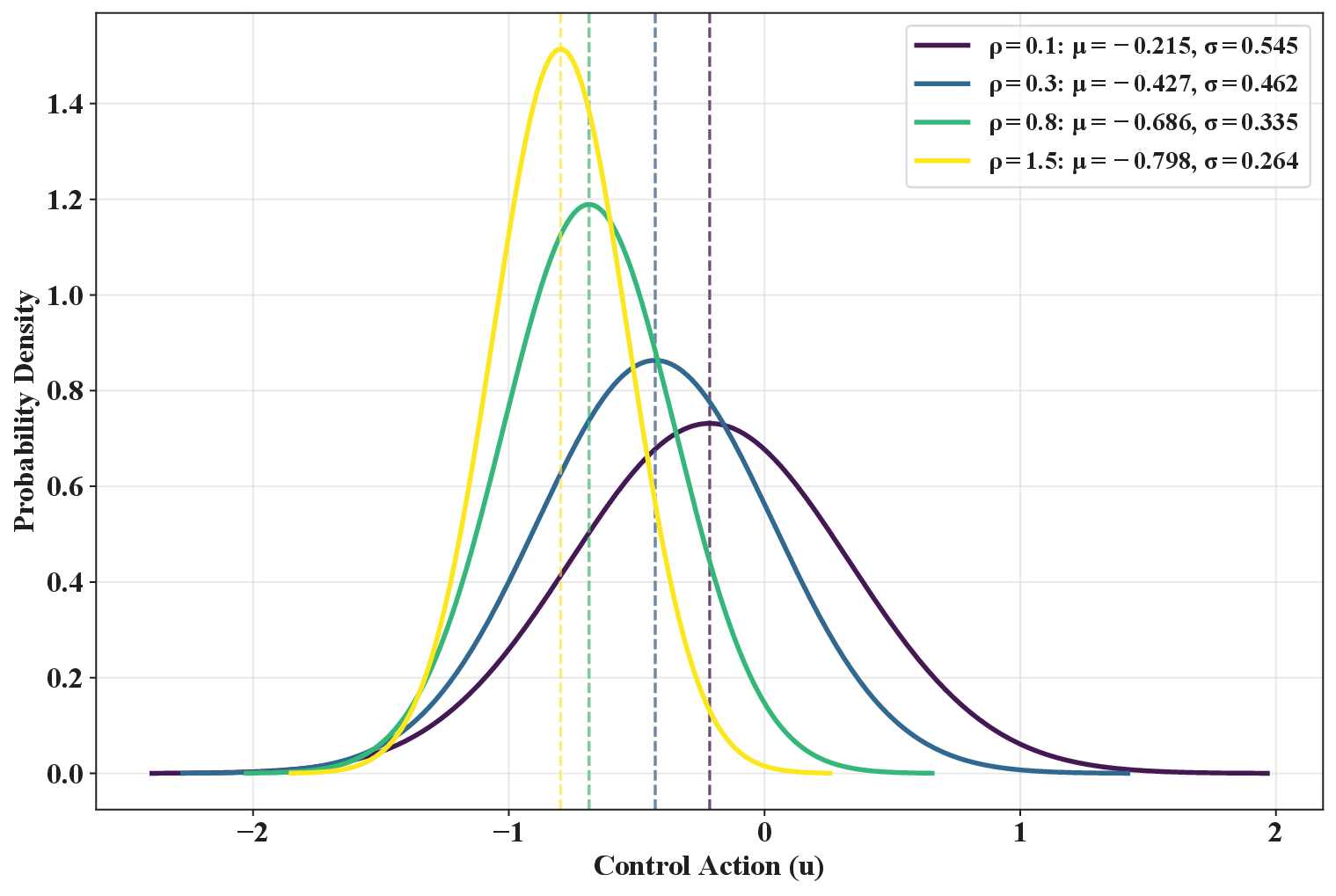} % 单栏适配：宽度改为0.8\linewidth
		\caption{The impact of $\rho$ on the optimal policy}
		\label{fig:3}
	\end{figure}

	In Fig.~\ref{fig:3}, with $\lambda$ fixed at 0.6 and $\overline{\sigma}$ fixed at 1.0, we analyze the impact of the discount rate $\rho$ on the optimal policy for the indoor temperature regulation system. The selected $\rho$ values (0.1, 0.3, 0.8, and 1.5) all satisfy both the solvability condition $\rho > \alpha$ and the concavity requirement $k_2 < 0$ of the value function.
	
	As demonstrated in Fig.~\ref{fig:3}, the discount rate $\rho$ effectively modulates the exploration variance of the optimal policy, consistent with the theoretical framework. Higher values of $\rho$ (e.g., $\rho = 1.5$), indicating greater myopia, yield more concentrated distributions with reduced variance ($\sigma^2 = 0.07$), reflecting diminished exploration incentives. Conversely, lower values of $\rho$ (e.g., $\rho = 0.1$), representing far-sighted optimization, produce broader distributions with increased variance ($\sigma^2 = 0.28$), facilitating more extensive exploration. This modulation is channeled through the coefficient $k_2$ in (18), as evidenced by the systematic variation of both the mean $\mu$ and variance $\sigma^2$ with $\rho$.
	
	\subsection{Convergence of Optimal Policy and Exploratory Value Function}
\begin{figure}[!h]
		\centering
		\includegraphics[width=0.5\linewidth]{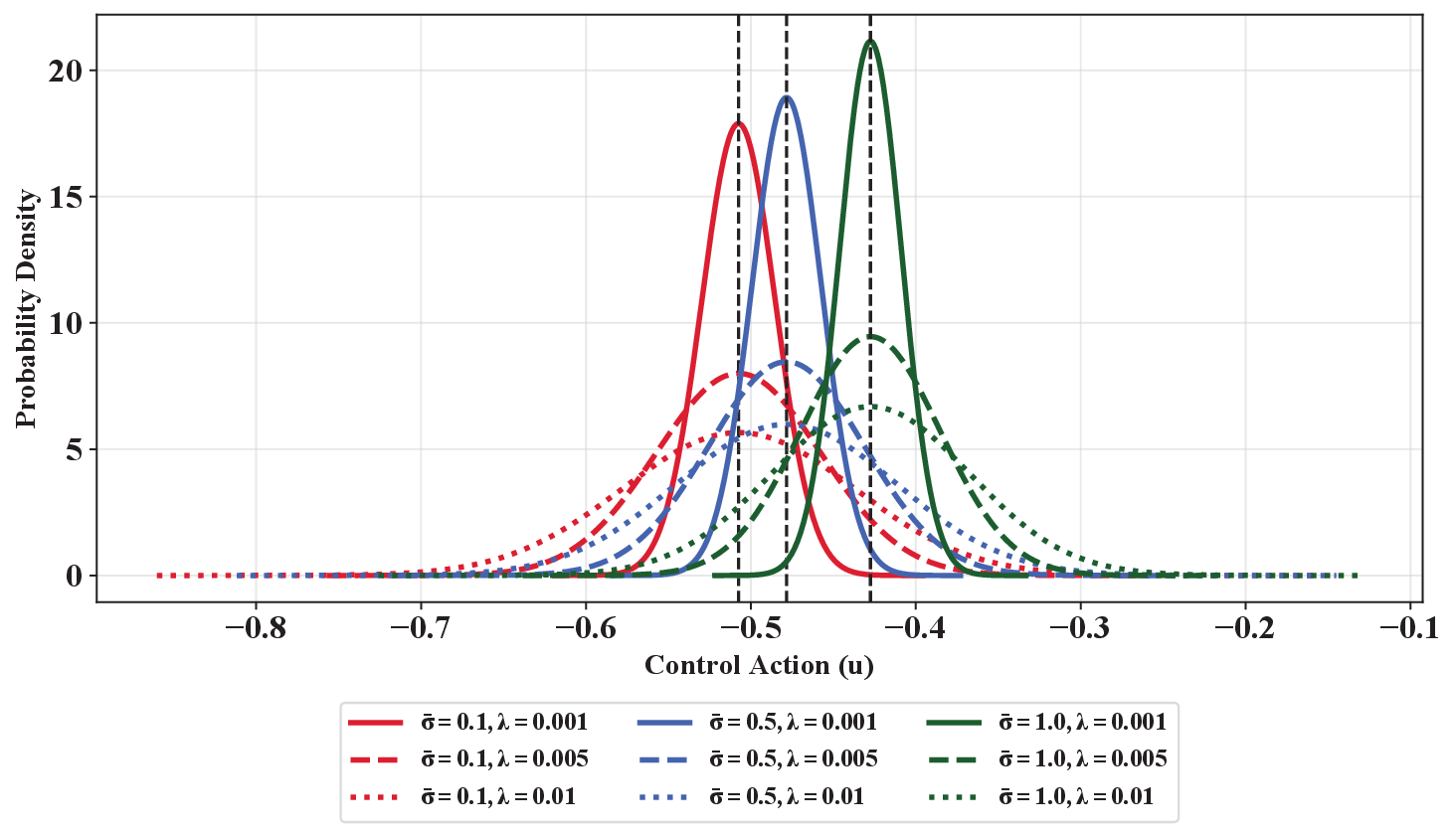} % 单栏适配：宽度改为0.8\linewidth
		\caption{Convergence of the optimal policy}
		\label{fig:4}
	\end{figure}

	In Fig.~\ref{fig:4}, we examine the impact of $\lambda$ (set to 0.001, 0.005, 0.01, and 0) on optimal policy convergence under varying $\overline{\sigma}$ values (0.1, 0.5, 1.0). Fig.~\ref{fig:4} investigates the asymptotic behavior of the exploratory value function as $\lambda \to 0^+$ (through a sequence of diminishing $\lambda \neq 0$ values), testing its convergence to the non-exploratory value function across the same $\overline{\sigma}$ range.
	
	Fig. \ref{fig:4} displays Gaussian distribution curves corresponding to different $\lambda$ values under varying $\overline{\sigma}$. A clear trend emerges: as $\lambda$ decreases, the distributions become progressively narrower (indicating reduced variance) with increasing kurtosis, ultimately converging to the mean. This demonstrates that the Gaussian distribution characterizing the optimal policy approaches a single-point distribution when $\lambda \rightarrow 0$.
	
	Based on Proposition 2, we establish the value functions for both the exploratory and non-exploratory problems. As shown in the figure, the dashed lines represent the value functions of the non-exploratory problems ($V^{ne}(x)$), while the solid lines correspond to the value functions of the exploratory problems ($V(x)$), with different colors indicating different Knightian uncertainty levels ($\overline{\sigma} = 0.1$, $0.5$, and $1.0$). The results demonstrate that as $\lambda \rightarrow 0$, the value functions of the exploratory problems converge to their non-exploratory counterparts across all Knightian uncertainty levels.
	
	The simulation results in Figs. \ref{fig:3} and \ref{fig:4} prove the theoretical result of Proposition 4 in the paper, i.e., when $\lambda \rightarrow 0$, the Gaussian distribution obeyed by the optimal policy converges to a single-pointed distribution and the value function of the exploratory problem converges to that of the non-exploratory problem, i.e., the exploratory LQ problem under Knightian uncertainty will converge to its non-exploratory counterpart.
	
	% ========== 附录内容结束 ==========

\end{document}